\documentclass[12pt]{article}
\usepackage{amssymb, amsmath, amsfonts}
\usepackage{mathrsfs}
\usepackage{dsfont}
\newcommand{\p}{{\sf P}}
\newcommand{\e}{{\sf E}}
\newcommand{\ind}{{\mathbb{I}}}

\newcommand{\var}{{\sf var}}
\newcommand{\cov}{{\sf cov}}

\newtheorem{Thm}{Theorem}

\newtheorem{remk}{Remark}

\newtheorem{cor}{Corollary}

\newtheorem{Lem}{Lemma}

\newcommand\norm[1]{\left\lVert#1\right\rVert}

\begin{document}

\begin{center}
{\LARGE{\bf Statistical estimation of the Shannon entropy }}
\end{center}

\begin{center}
{\large {\bf Alexander Bulinski\footnote{E-mail: bulinski@yandex.ru}, Denis Dimitrov\footnote{E-mail: den.dimitrov@gmail.com}}}
\end{center}
\begin{center}
{\it Steklov Mathematical Institute of Russian Academy of Sciences;\\ Dept. of Mathematics and Mechanics, Lomonosov Moscow State University,\\ Moscow 119234, Russia}
\end{center}
\vskip0.5cm
{\small
{\bf Abstract} The behavior of the Kozachenko - Leonenko estimates for the (differential) Shannon entropy is studied when the number of i.i.d. vector-valued observations tends to infinity. The asymptotic unbiasedness and $L^2$-consistency of the estimates are established. The conditions employed involve the analogues of the Hardy - Littlewood maximal function. It is shown that the results are valid in particular for the entropy estimation of any nondegenerate Gaussian vector.}

\vskip0.2cm
\noindent
{\small {\bf Key words}  Shannon differential entropy; Kozachenko - Leonenko estimates;  Hardy - Littlewood maximal function analogues; asymptotic unbiasedness and $L^2$-consistency; Gaussian vectors}
\vskip0.2cm
\noindent
{\bf MR (2010) Subject Classification} 60F25, 62G20, 62H12

\section{Introduction}
The notion of entropy belongs to the principle ones in Physics and Mathematics.
R.Clausi\-us is considered as the father of the entropy concept.
The important contributions to the development of this concept
were made by L.Boltzmann, J.Gibbs and M.Planck.  Mathematicians also were preoccupied with the entropy.  The works by C.Shannon, A.N.Kolmogorov, Ya.G.Sinai, A.R\'enyi, A.S.Holevo, T.Tsallis are worth mentioning in this regard. On the history of different paths to entropy see, e.g.,  \cite{Balibrea}, \cite{Benguigui}.
Thus there are various definitions of entropy. We recall that proposed by C.Shannon
in \cite{Shannon}. Namely, if a random variable $\xi$
takes values in
a finite or numerable set $S$, then the Shannon entropy of $\xi$ (or
of $\xi$ distribution)
is given by the formula
$$
H:=- \sum_{x\in S}p_x\log p_x
$$
where $p_x:=\p(\xi=x)$, $x\in S$. Here the logarithm base is $2$ but we will employ the base $e$ since a constant factor is not essential below. As usual we set $0 \log 0:=0$. Now, for a random vector $\xi$ taking values in $\mathbb{R}^d$ and having a density $f$ w.r.t. the Lebesgue measure $\mu$, one can introduce in a similar way the entropy (also called differential entropy)
\begin{equation}\label{1}
H(f):= -\int_{\mathbb{R}^d}f(x)\log f(x) \mu(dx).
\end{equation}
For a density version $f$
denote the support $S(f):=\{x\in \mathbb{R}^d: f(x)>0\}$. Clearly, the integral in \eqref{1}
is taken over $S(f)$.
To simplify notation we will
write $dx$ instead of $\mu(dx)$.

We mention in passing the R\'enyi entropy $H_{\alpha}$ and the Tsallis entropy $S_q$  depending on certain parameters $\alpha$ and $q$, respectively (when $\alpha\to 1$ and $q\to 1$  the Shannon entropy arises).
On these and other divergence measures see, e.g., \cite{Gorban}.

Statistical estimates of $H(f)$ constructed by means of i.i.d. observations $X_1,\ldots,X_N$ having the same law as $\xi$ are very important.
They permit to estimate some other characteristics of distributions, e.g., the mutual information for two random vectors.
Such estimates are widely used in machine learning, they are
employed in entropy-based goodness-of-fit tests,
they are applied in the feature selection theory and in the detection of  texture inhomogeneities (see, e.g., \cite{Alonso-Ruiz}, \cite{Bierlant}, \cite{Leonenko}, \cite{Peng}).
We leave apart many other domains where an entropy (entropies) estimates are quite useful.

There exist several approaches to $H(f)$ estimation. We are interested in nonparametric models.
Note that $H(f)=\e (-\log f(\xi))$. The ``plug-in'' method leads to the estimates
$\widehat{H}_N:= -\frac{1}{N}\sum_{k=1}^N \log \widehat{f}(X_k)$ where $\widehat{f}(x)$ is a
(kernel) estimate of a density $f$ at a point $x\in \mathbb{R}^d$.
In this regard we refer to the review \cite{Bierlant}.
Considering $H(f)$ estimation problem E.G.Miller writes in \cite{Miller} that
plug-in estimates work well in low dimensions  and for densities with known parametric form,
and the difficult problem of density estimation makes them impractical for small sample sizes in
higher dimensions.  He introduced the $m$-Voronoi and $m$-Delaunay estimators for $H(f)$,
however the theorems concerning these estimates behavior were not provided.
The method of a recursive rectilinear partitioning for $H(f)$ estimation is proposed in
\cite{Stowell}. The author compares his method with Miller's one
and discusses the complexity of various partitioning schemes.
J.Ma and Z.Sun \cite{Ma} analyzed the copula entropy estimation  to get the mutual information estimate.
The Gaussian copula in the entropy estimation is also used  in \cite{Charzynska}.
A nonparanormal information estimation is considered in \cite{Singh}.
An approach involving the nearest neighbor statistics to estimate an entropy was proposed by L.F.Kozachenko and N.N.Leonenko in a well-known paper \cite{Kozachenko} and developed by N.N.Leonenko with coauthors in a series of papers (see, e.g., \cite{Leonenko}).
These estimates can be viewed as the basis for the widely applied mutual information estimates
introduced by A.Kraskov et al. \cite{Kraskov}.
During the last two decades many authors
employed and analyzed
the Kozachenko - Leonenko estimates. Unfortunately, in a number of papers the proofs of asymptotic properties of $H(f)$ estimates  were not correct as it was pointed out, e.g., in \cite{Pal}.

The main goal of our research is to provide, under wide conditions, the proof of the asymptotic unbiasedness of $H(f)$ estimates and, moreover, to establish their $L^2$-consistency. Note that paper \cite{Pal} is devoted to proving the strong consistency of the specified estimates of the R\'enyi entropy $H_{\alpha}$ for $\alpha \in (0,1)$  when $\xi$ has a density $f$ with bounded support. If  $f$ is a Lipschitz function, then for any $\delta \in (0,1)$ with probability at least $1-\delta$ the convergence rate (depending on $d$ and $\delta$) of these estimates to $H_{\alpha}$ is provided in \cite{Pal} as well. The authors of \cite{Pal} indicated that their theory would need
significant changes
to comprise the Shannon entropy estimation.
For certain estimates, the convergence rate to $H(f)$ and the convergence in distribution to the normal law under appropriate normalization are established, e.g., in \cite{Delattre}, \cite{Laurent}. The proofs exploit assumption that a density $f$ is rather smooth. We do not impose such requirement.
To complete the brief introduction we mention that, according to \cite{Sricharan}, the ensemble methods using the combinations of weighted statistics can improve the convergence rate for
initially constructed statistics. Note also that it is interesting to consider estimates of
$H(f)$ involving $k$-NN statistics (see, e.g., \cite{Charzynska} and \cite{Singh_1} where the authors use 
conditions on $f$ different from those employed here).

The paper is organized as follows. In Section 2 we introduce the estimates $H_N$ of
$H(f)$ and formulate two main results. Section 3 is devoted to the proof of
$H_N$ asymptotical unbiasedness. The $L^2$-consistency of $H_N$ is demonstrated in Section 4.
The proofs of some auxiliary results are given in Appendix.

\section{Main results}
Let $X_1,\ldots,X_N$ be i.i.d. random vectors  having the same law as a vector $\xi$ with values in $\mathbb{R}^d$. All random vectors (variables) under consideration are defined on a complete probability space $(\Omega,\mathcal{F},\p)$. Assume that $\xi$ (i.e. the distribution $\p_{\xi}$ of $\xi$) has a density $f$ w.r.t. the Lebesgue measure $\mu$ in $\mathbb{R}^d$. For each $i=1,\ldots,N$, set $\rho_i:=\min\{\rho(X_i,X_j): j \in \{1,\ldots,N\}\setminus \{i\}\}$, where $\rho(x,y)$ is the Euclidean distance between $x,y\in \mathbb{R}^d$. In other words $\rho_i$ is the distance from $X_i$ to its nearest neighbor in the sample $\{X_1,\ldots,X_N\}\setminus \{X_i\}$. Further on we consider $N\geq 2$. Introduce $\overline{\rho}:= (\rho_1\cdot \ldots \cdot\rho_N)^{1/N}$. Recall that the  Kozachenko-Leonenko estimate of an entropy $H$ is provided by the formula
\begin{equation}\label{KLestimation}
H_N:=d \log{\overline{\rho}} + \log{V_d} + \gamma + \log{(N-1)}
\end{equation}
where $\gamma := -\int_{(0,\infty)} e^{-t} \log{t}\, dt \; \approx 0.5772$ and $V_d := \frac{{\pi}^{\frac{d}{2}}}{\Gamma(\frac{d}{2}+1)}$ are the Euler constant
and   the volume (i.e. the Lebesgue measure) of a unit ball in
$\mathbb{R}^d$, respectively. We can write $H_N= \log(\overline{\rho}^d V_d \widetilde{\gamma}(N-1))$ with
$\widetilde{\gamma}:=\exp\{\gamma\}$. Let $B(x,r) = \{y \in \mathbb{R}^d: \rho(x,y) \leq r\}$ be the ball of a radius $r\geq 0$ with a center $x\in \mathbb{R}^d$.
Clearly, its volume $|B(x,r)| := \mu(B(x,r)) = r^d V_d$. Set
\begin{equation}\label{G}
G(t) :=
\left\{
\begin{aligned}
& 0, \;\, \;\; \; 0 \leq t < 1, \\
& t \log{t}, \;\, \;\; \; t \geq 1.
\end{aligned}
\right.
\end{equation}
For a probability density $f$ in $\mathbb{R}^d$,
$x\in \mathbb{R}^d$, $r>0$ and $R>0$, introduce the functions (or functionals depending on parameters)
\begin{equation}\label{I}
I_f(x,r):= \frac{\int_{B(x,r)} f(y) \, dy}{r^d V_d},
\end{equation}
\begin{equation}\label{mM}
M_f(x, R) := \sup_{r \in (0,R]} I_f(x,r),\;\;m_f(x, R) := \inf_{r \in (0,R]} I_f(x,r).
\end{equation}
We use the following elementary result.
\begin{Lem}\label{l6a} For a probability density $f$ in $\mathbb{R}^d$, the function $I_f(x,r)$ defined in \eqref{I} is continuous in $(x,r)\in \mathbb{R}^d\times (0,\infty)$.
\end{Lem}
The proof is provided in Appendix.
Thus, invoking, e.g., Theorem 15.84  \cite{Yeh}, we can claim that, for each $R>0$, the functions
$m_f(\cdot,R)$ and $M_f(\cdot,R)$ are upper semicontinuous and lower semicontinuous, respectively.
Hence, in view of Proposition 15.82 \cite{Yeh} these nonnegative functions are Borel measurable. Clearly, for each $x\in \mathbb{R}^d$, $m_f(x,\cdot)$ is nonincreasing and $M_f(x,\cdot)$ is nondecreasing.
Note in passing that changing $\sup_{r \in (0,R]}$ by $\sup_{r \in (0,\infty)}$
in the definition of $M_f(x,R)$ leads to the celebrated Hardy - Littlewood maximal
function $M_f(x)$ widely used in the harmonic analysis.
Some properties
of a function $\int_{B(x,r)} f(y) \, dy$ are considered, e.g., in \cite{Evans}.

For a probability density $f$ in $\mathbb{R}^d$, positive $\varepsilon_i,R_j$, where $i=0,1,2$ and $j=1,2$, we define the following functionals
with values in $[0,\infty]$
\begin{equation}\label{p1}
K_f(\varepsilon_0):=\int_{\mathbb{R}^d} {\Big(\int_{\mathbb{R}^d} {G\big(|\log{\rho}(x,y)|\big)} f(y) \,dy \Big)}^{1+\varepsilon_0} f(x) \, dx,
\end{equation}
\vspace{-0.2cm}
\begin{equation}\label{p2}
Q_f(\varepsilon_1, R_1):=\int_{\mathbb{R}^d} M_f^{\varepsilon_1}(x, R_1) f(x)\,dx,
\end{equation}
\vspace{-0.2cm}
\begin{equation}\label{p3}
T_f(\varepsilon_2,R_2):=\int_{\mathbb{R}^d} m_f^{-\varepsilon_2}(x, R_2) f(x) \, dx.
\end{equation}

\begin{remk}\label{rem1}
\normalfont{One has to write $\int_{\mathbb{R}^d\setminus \{x\}} {G\big(|\log{\rho}(x,y)|\big)} f(y) \,dy$ in \eqref{p1} as $\rho(x,x)=0$.  However, we can
formally set $\log 0:=-\infty$ and $G(\infty):=\infty$ to keep formula \eqref{p1} since
$\p_{\xi}(\{x\})=0$ for any $x\in \mathbb{R}^d$. We also suppose that
$1/0:=\infty$ (consequently $m_f^{-\varepsilon_2}(x, R_2):=\infty$ when $m_f(x,R_2)=0$).
For each version of a function $f\in L^1(\mathbb{R}^d)$, we can write in \eqref{p1}, \eqref{p2}, \eqref{p3} the integral over the support $S(f):=\{x\in \mathbb{R}^d: f(x)>0\}$  instead of integrating over $\mathbb{R}^d$ (evidently, the results do not depend on the choice of $f$ version).}
\end{remk}

\vskip0.2cm
\begin{Thm}\label{th1}
Assume that, for some positive $\varepsilon_i,R_j$, where $i=0,1,2$ and $j=1,2$, the functionals
appearing in \eqref{p1}, \eqref{p2}, \eqref{p3} are finite, so $K_f(\varepsilon_0) < \infty$, $Q_f(\varepsilon_1, R_1) < \infty$, $T_f(\varepsilon_2, R_2) < \infty$.
Then the estimates $H_N$ are asymptotically unbiased, i.e.
\begin{equation}\label{main1}
\lim\limits_{\scriptscriptstyle{N \rightarrow \infty}} \e H_N = H.
\end{equation}
\end{Thm}

\begin{remk}\label{rem2}
\normalfont{It is useful to note that if $Q_f(\varepsilon_1, R_1) < \infty$ and $T_f(\varepsilon_2, R_2) < \infty$ for some positive $\varepsilon_1, \, \varepsilon_2, R_1, R_2$ then $\int_{\mathbb{R}^d} |\log{f(x)}|f(x) \, dx < \infty$.
Indeed, definition \eqref{mM}
and the Lebesgue differentiation theorem (see, e.g., Theorem 25.17 \cite{Yeh})
yield that $m_f(x, R_2) \leq f(x) \leq M_f(x, R_1)$ for $\mu$-almost all $x \in \mathbb{R}^d$.
Evidently, $\log z \leq \frac{1}{\varepsilon}z^{\varepsilon}$ for any $z\geq 1$ and each $\varepsilon >0$. Consequently,
\begin{gather*}\int_{\mathbb{R}^d} |\log{f(x)}| f(x) \, dx = \int_{f(x) \geq 1} \log{f(x)} f(x) \, dx + \int_{f(x) < 1} (-\log{f(x)}) f(x) \, dx \\ \leq  \frac{1}{\varepsilon_1} Q_f(\varepsilon_1, R_1) + \frac{1}{\varepsilon_2} T_f(\varepsilon_2, R_2) < \infty.
\end{gather*}
So, the finiteness of integrals \eqref{p2} and \eqref{p3} implies $|H(f)| < \infty$.

}
\end{remk}

We formulate a simple but useful result.

\begin{Lem}\label{lemma1} Let $f$ be a probability density in $\mathbb{R}^d$.
Then the following statements are valid.

$1)$ If $K_f(\varepsilon_0)<\infty$ for some $\varepsilon_0>0$ then $K_f(\varepsilon)<\infty$ for any $\varepsilon \in (0,\varepsilon_0]$.

$2)$ If $Q_f(\varepsilon_1, R_1)<\infty$ for some $\varepsilon_1 > 0$ and $R_1>0$ then $Q_f(\varepsilon, R)<\infty$ for any $\varepsilon \in (0,\varepsilon_1]$ and each $R>0$.

$3)$ If $T_f(\varepsilon_2,R_2)<\infty$ for some $\varepsilon_2>0$ and $R_2>0$ then $T_f(\varepsilon,R)<\infty$ for any $\varepsilon \in (0,\varepsilon_2]$ and each $R>0$.
\end{Lem}

The proof of this Lemma  is provided in Appendix. In view of Lemma \ref{lemma1} one can recast
Theorem \ref{th1} as follows.

\begin{cor}\label{cor1} Let $f$ be a probability density in $\mathbb{R}^d$ such that,
for some $\varepsilon >0$,
$K_f(\varepsilon)$, $Q_f(\varepsilon)$ and $T_f(\varepsilon)$ are finite where
$Q_f(\varepsilon):=Q_f(\varepsilon, \varepsilon)$, $T_f(\varepsilon):=T_f(\varepsilon,\varepsilon)$. Then \eqref{main1} holds.
\end{cor}

Let us consider the following conditions.
\vskip0.3cm
\noindent
(A) For some $p>1$,
\begin{equation}\label{p}
\int_{\mathbb{R}^d}\int_{\mathbb{R}^d}|\log \rho(x,y)|^p f(x)f(y)\,dxdy<\infty.
\end{equation}
(B) There exists a version of density $f$ such that, for some $M>0$,
\begin{equation*}
f(x)\leq M<\infty,\;\;x\in \mathbb{R}^d.
\end{equation*}
(C1) There exists a version of density $f$ such that, for some $m>0$,
\begin{equation*}
f(x)\geq m >0,\;\;x\in S(f).
\end{equation*}

\begin{cor}\label{cor2}
Any assumption of Theorem \ref{th1} concerning the finiteness of integrals \eqref{p1}, \eqref{p2}, \eqref{p3} can be replaced by conditions {\normalfont (A), (B), (C1)}, respectively, and then \eqref{main1} will be true as well. Moreover, if {\normalfont (B)} and {\normalfont (C1)} are satisfied then \eqref{main1} holds whenever $f$ has a bounded support.
\end{cor}

The proof of this Corollary is provided in Section 3.

Along with \eqref{p1}, for a probability density $f$ in $\mathbb{R}^d$ and positive $\varepsilon_0$, we define the following functional with values in $[0, \infty]$
\begin{equation}\label{p12}
K_{f,2}(\varepsilon_0):=\int_{\mathbb{R}^d} {\Big(\int_{\mathbb{R}^d} {G\big(\log^2{\rho}(x,y)\big)} f(y) \,dy \Big)}^{1+\varepsilon_0} f(x) \, dx.
\end{equation}
Note that the statement $1)$ of Lemma \ref{lemma1} is true for $K_{f,2}$ as well, so if $K_{f,2}(\varepsilon_0) < \infty$ then $K_{f,2}(\varepsilon) < \infty$ for any $\varepsilon \in (0, \varepsilon_0]$ (see the proof of  Lemma \ref{lemma1}).
Let us formulate the conditions that guarantee $L^2$-consistency of \eqref{KLestimation}.
\begin{Thm}\label{th_main2}
Assume that, for some positive $\varepsilon_i,R_j$ where $i=0,1,2$, $j=1,2$, the functionals
appearing in \eqref{p12}, \eqref{p2}, \eqref{p3} are finite, so $K_{f,2}(\varepsilon_0) < \infty$, $Q_f(\varepsilon_1, R_1) < \infty$, $T_f(\varepsilon_2, R_2) < \infty$.
Then the estimates $H_N$ are $L^2$-consistent, i.e.
\begin{equation}\label{main2}
\e(H_N-H)^2 \to 0,\;\;N\to \infty.
\end{equation}

\end{Thm}

\begin{cor}\label{cor3}
Any assumption of Theorem $\ref{th_main2}$ concerning the finiteness of integrals \eqref{p12}, \eqref{p2} and \eqref{p3} can be replaced, respectively, by the following ones: \eqref{p} is valid for some $p > 2$, {\normalfont (B) and (C1)}. Then \eqref{main2} will be true. Moreover, if {\normalfont (B)} and {\normalfont (C1)} are satisfied then \eqref{main2} holds whenever $f$ has a bounded support.
\end{cor}

Now instead of (C1) we consider the following condition.
\vskip0.3cm
\noindent
(C2) For a fixed $R>0$, there exists a constant $c>0$ and a version of a density $f$ such that
\begin{equation}\label{reg}
m_f(x,R)\geq c f(x),\;\;x\in \mathbb{R}^d.
\end{equation}

Note that D.Evans  considered the ``positive density condition''
in Definition 2.1 of \cite{Evans}
meaning that there exist constants $\beta >1$ and $\delta >0$ such that $\frac{r^d}{\beta}\leq \int_{B(x,r)}f(y)dy \leq \beta r^d $ for all $0\leq r\leq \delta$ and $x\in \mathbb{R}^d$.
Consequently $m_f(x,\delta)\geq \frac{1}{\beta V_d}$, $x\in \mathbb{R}^d$.
It was proved in \cite{Evans_1} that if $f$ is smooth and its support is a compact convex body in $\mathbb{R}^d$ then the mentioned inequalities from Definition 2.1 of \cite{Evans} hold.

\begin{remk}\label{rem2a}
\normalfont{
If, for some positive $\varepsilon$, $R$ and $c$, condition (C2) is true and
\begin{equation}\label{densi}
\int_{\mathbb{R}^d}f(x)^{1-\varepsilon}dx <\infty
\end{equation}
then obviously $T_f(\varepsilon,R)<\infty$.
Thus in Theorems \ref{th1} and \ref{th_main2} we can employ (C2) and \eqref{densi} instead of the assumption
concerning $T_f(\varepsilon,R)$. To illustrate this observation
we provide the following result for a density with unbounded support.}
\end{remk}

\begin{cor}\label{cor4}
Let $\xi$ be a Gaussian random vector in $\mathbb{R}^d$ with $\e \xi = \nu$ and a nondegenerate covariance matrix $\Sigma$ $($i.e. $\xi$ has a density$)$. Then relations \eqref{main1} and \eqref{main2} hold where
$H= \frac{1}{2}\log {\sf det}(2\pi e \Sigma)$.
\end{cor}

The proofs of Corollaries \ref{cor3} and \ref{cor4} are given in Section 4.

\section{Proof of Theorem 1}
\vskip0.3cm

According to the Lebesgue differentiation theorem  if $f\in L^1(\mathbb{R}^d)$ then, for $\mu$-almost all  $x\in \mathbb{R}^d$,
the following relation holds
\begin{equation}\label{2}
\lim_{r\to 0+}\frac{1}{|B(x,r)|}\int_{B(x,r)}|f(y)-f(x)|\,dy=0.
\end{equation}
Let $\Lambda(f)$ stand for a set of all the Lebesgue points of a function $f$,
i.e. for $x\in \mathbb{R}^d$ satisfying \eqref{2}. Clearly $\Lambda(f)$ depends on the chosen
 version of $f$ belonging to the class of equivalent functions from $L^1(\mathbb{R}^d)$.
For each version of $f$ we have $\mu(\mathbb{R}^d\setminus \Lambda(f))=0$.

We can rewrite the estimate $H_N$ as follows
$$
H_N = \frac{1}{N} \sum\limits_{i=1}^N \zeta_i(N),  \;\; \zeta_i(N) = \log{\big(\rho_i^d V_d \widetilde\gamma (N-1)\big)},\;\;i=1,\ldots,N.
$$
Recall that $\rho_i\!=\!\min\{\rho(X_i,X_j)|\!:\! j \!\in\! \{1,\ldots,N\}\setminus \{i\}\!\}$. The random variables
$\zeta_1(N),\ldots,\zeta_N(N)$ are identically distributed since
$X_1,\ldots,X_N$ are i.i.d. random vectors, and therefore
to prove Theorem \ref{th1} we will show that
$\e |\zeta_1(N)|<\infty$ for all $N$ large enough and
\begin{equation}\label{mr}
\e \zeta_1(N)\to H,\;\;N\to \infty.
\end{equation}
Note that if $V$ is a nonnegative (a.s.) random variable (hence $\e V \leq \infty$) and $X$ is an arbitrary random vector with values in $\mathbb{R}^d$ then there exists $\e (V|X)$. Moreover, $\e (V|X) = \psi(X)$ where a measurable function $\psi$ takes values in $[0,\infty]$ and
\begin{equation}\label{ce}
\e V=\int_{\mathbb{R}^d} \e (V|X=x)\p_X(dx),
\end{equation}
i.e. $\e (V|X=x):=\psi(x)$, $x\in \mathbb{R}^d$. Formula \eqref{ce} means that
simultaneously both sides are finite or infinite and coincide (to verify \eqref{ce} one can use the usual formula
for a random variables $V_k := V\ind\{k\leq V<k+1\}$ with finite expectation and the representation $V=\sum_{k=0}^{\infty}V_k$ for nonnegative (a.s) random variable $V$ with $\e V \leq \infty$). Let $F(u,\omega)$ be a regular conditional distribution function of $V$ given $X$ where $u \in [0,\infty)$ and $\omega \in \Omega$, i.e.  $F(u,\omega)$ is the specified version of $\p(V\leq u|X)(\omega)$, see, e.g.,  Theorem 4, Ch. 2, Sect. 7 in \cite{Shiryaev}.
We will consider only such version since it provides the possibility to evaluate the conditional expectation according to the extension of Theorem 3, Ch. 2, Sect.~7 in \cite{Shiryaev}.
Namely, if $h$ is a measurable function such that $h:\mathbb{R}\to [0,\infty)$ then, for $\p_{X}$-almost all $x\in \mathbb{R}^d$, one has without assumption $\e h(V)<\infty$ that
\begin{equation}\label{ce1}
\e (h(V)|X=x)= \int_{[0,\infty)}h(u)dF(u,x).
\end{equation}
It means that both sides of \eqref{ce1} are finite or infinite simultaneously and coincide
(if $\e h(V)=\infty$, to establish \eqref{ce1} we consider $h_n(u):= h(u)\ind\{0\leq u\leq n\}$, $n\in \mathbb{N}$, and use the monotone convergence theorem).
By means of \eqref{ce} and \eqref{ce1} one can prove that $\e |\zeta_1(N)|<\infty$ for all $N$ large enough and \eqref{mr} holds. For this purpose we take $V=e^{\zeta_1(N)}$,
$X=X_1$ and $h(u)= |\log u|$, $u>0$ (we take $h(u)=\log^2 u$ in the proof of Theorem \ref{th_main2}). Writing $\log U$ for a positive a.s. random variable $U$
we set $\log U(\omega):=-\infty$ when $U(\omega)=0$, and as usual we stipulate that
$\int_A g(y) Q(dy)=0$ whenever $g(y)= -\infty$ (or $+\infty$) for $y\in A$ and $Q(A)=0$. To reduce the volume of the paper we
only consider below the evaluation of $\e \zeta_1(N)$ as all steps of the proof are the same when
treating $\e|\zeta_1(N)|$.

We divide the proof of Theorem \ref{th1} into four steps. Steps 1-3  are devoted to the demonstration of relation
\begin{equation}\label{convc}
\e(\zeta_1(N)|X_1=x)\to -\log f(x), \;\;x\in A\subset S(f),\;\;N\to \infty,
\end{equation}
where $A$ depends on an $f$ version and $\p_{\xi}(S(f)\setminus A)=0$.
Step 4 contains the proof of the desired statement \eqref{mr}.

\vskip0.3cm
{\it Step 1}.
For $x\in \mathbb{R}^d$ and $u>0$ we study the asymptotic behavior (as $N\to \infty$) of the following function
\begin{align}\label{eq1a}
\begin{gathered}
    F_{N,x}(u) := \p\big(e^{\zeta_1(N)} \leq u | X_1 = x\big) = \p\big(\rho_1^d(N) V_d \widetilde\gamma (N-1) \leq u | X_1 = x\big)\\
    = \p\left(\min_{j =2,\ldots,N} \rho(X_1, X_j) \leq \textstyle{{\left(\frac{u}{V_d \widetilde\gamma (N-1)}\right)}^{\frac{1}{d}}} \big| X_1 = x\right)
    \\= \p\big( \min_{j =2,\ldots,N} \rho(x, X_j) \leq r_N(u)  \big)= \p(\xi_{N,x}\leq u)
\end{gathered}
\end{align}
where
\vspace{-0.4cm}
\begin{equation}\label{rN}
r_N(u):= {\left(\frac{u}{V_d \widetilde\gamma (N-1)}\right)}^{\frac{1}{d}},
\end{equation}
\begin{equation}\label{eq3a}
\xi_{N,x}:=(N-1)V_d\widetilde{\gamma}\min_{j =2,\ldots,N} \rho^d(x, X_j),
\end{equation}
and we have employed in \eqref{eq1a} the independence of random vectors $X_1,\ldots,X_N$.
Therefore
\begin{align}\label{a1}
\begin{gathered}
F_{N,x}(u) = 1 - P\big( \min_{j =2,\ldots,N} \rho(x, X_j) > r_N(u)  \big)
\\
= 1 - {\Big(1 - P\big(\xi \in B(x, r_N(u))\big)\Big)}^{N-1}
= 1 - {\Big(1 -  \int\limits_{B(x,r_N(u))} f(y) \, dy  \Big)}^{N-1}
\end{gathered}
\end{align}
because $X_1,\ldots,X_N$ are independent copies of a vector $\xi$. Formula \eqref{a1} shows that $F_{N,x}(u)$ is
the regular conditional distribution function of $e^{\zeta_1(N)}$ given $X_1=x$.

Note that, for each $u>0$, $r_N(u)  \rightarrow 0$ as $N \rightarrow \infty$ and $|B(x,r_N(u))| = V_d {\big(r_N(u)\big)}^d = \frac{u}{\widetilde\gamma (N-1)}$. Hence in view of  \eqref{2}, for any fixed $x\in \Lambda(f)$ and $u> 0$,
$$
\frac{\widetilde\gamma (N-1)}{u} \int\limits_{B(x,r_N(u))} f(y) \, dy = f(x) + \alpha_N{(x,u)}
$$
where $\alpha_N{(x,u)} = \bar{\bar{o}}(1), \; N \rightarrow \infty$. Thus
\eqref{a1} implies that, for $u>0$ and $x\in \Lambda(f)\cap S(f)$, we get (as $f(x)>0$ for $x\in S(f)$)
\begin{align}\label{conv}
\begin{gathered}
    \lim_{\scriptscriptstyle{N \rightarrow \infty}} F_{N,x}(u) =  1 - \lim_{\scriptscriptstyle{N \rightarrow \infty}} {\bigg(1 - \frac{u}{\widetilde\gamma (N-1)} \Big(f(x) + \alpha_N{(x,u)}\Big) \bigg)}^{N-1} \\
    = 1 - e^{-\frac{f(x) u}{\widetilde\gamma}} =: F_{x}(u)=\p(\xi_x\leq u)
\end{gathered}
\end{align}
where $\xi_x \sim Exp\left(\frac{f(x)}{\widetilde{\gamma}}\right)$, $x\in S(f)$.
Relation \eqref{conv} means that
\begin{equation}\label{claw}
\xi_{N,x}\stackrel{law}\rightarrow \xi_x,\;\;x\in  \Lambda(f)\cap S(f),\;\;N\to \infty.
\end{equation}
We assume w.l.g. that, for all $x\in S(f)$, the random variables $\xi_x$ and $\{\xi_{N,x}\}_{N\geq 2}$
are defined on a probability space $(\Omega,\mathcal{F},\p)$ since in view of the Lomnicki - Ulam theorem (see, e.g. \cite{Kallenberg}, p. 93) one can consider
the independent copies of $X_1,X_2,\ldots$ and $\{\xi_x\}_{x\in S(f)}$ defined on a certain probability space.
The convergence in law of random variables is preserved under continuous mapping. Hence, for any
$x\in  \Lambda(f)\cap S(f)$, we come to the relation
\begin{equation}\label{b2}
\log \xi_{N,x}\stackrel{law}\rightarrow \log \xi_x,\;\;N\to \infty.
\end{equation}
We took into account that, for each $x\in \Lambda(f)\cap S(f)$, one has $\xi_x>0$ a.s. and
since $\xi$ has a density we infer that
$$
\p(\xi_{N,x}>0)
= (\p(\rho(x,\xi)>0))^{N-1}=(1-\p(\xi=x))^{N-1}=1,\;\;N\geq 2.
$$
More precisely, we can ignore zero values of nonnegative random variables (having zero values with probability zero)
when we take their logarithms.
\vskip0.3cm
{\it Step 2}. Now we will prove that, for $\mu$-almost every $x\in S(f)$,
\begin{equation}\label{a2}
\e \log \xi_{N,x}\to \e \log \xi_x,\;\;N\to \infty.
\end{equation}
Note that if $\eta\sim Exp(\lambda)$ where $\lambda >0$ then
$$
\e \log \eta = \int_{(0,\infty)}(\log u) \, \lambda e^{-\lambda u}\,du=
\int_{(0,\infty)}\log v e^{-v} dv - \log \lambda = -\log (\lambda \widetilde{\gamma}).
$$
Consequently, as $f(x)>0$ for $x\in S(f)$ and $\lambda=\frac{f(x)}{\widetilde{\gamma}}$, we get $\e \log \xi_{x}= - \log f(x)$.
It is easily seen that, for each $x\in \mathbb{R}^d$,
\begin{gather*}
    \e \log{\xi_{N,x}} = \int_{(0,\infty)} \log{u} \, dF_{N,x}(u) = \int_{(0,\infty)} \log{u} \, d\p(e^{\zeta_1(N)} \leq u|X_1 = x)  \\
    = \e(\log{e^{\zeta_1(N)}} | X_1 = x) = \e(\zeta_1(N) | X_1 = x).
\end{gather*}
Thus, for $x\in \Lambda(f) \cap S(f)$, the relation $\e(\zeta_1(N)|X_1=x)\to -\log f(x)$  holds if and only if $\e \log \xi_{N,x}\to \e \log \xi_x$ ($N\to \infty$).

According to Theorem 3.5 \cite{Billingsley} we  will establish \eqref{a2} if relation \eqref{b2} can be
accompanied, for $\mu$-almost all $x\in S(f)$, by the uniform integrability of a family $\{\log \xi_{N,x}\}_{N\geq N_0(x)}$.
Note that  $G$ introduced by \eqref{G} is an increasing function on $(0,\infty)$ and $\frac{G(t)}{t}\to \infty$ as $t\to \infty$. Therefore,
by the de la Valle Poussin theorem (see, e.g., Theorem 1.3.4 \cite{Borkar})
to guarantee, for $\mu$-almost every $x\in S(f)$, the uniform integrability of $\{\log \xi_{N,x}\}_{N\geq N_0(x)}$
it is sufficient to prove for such $x$, a positive $C_0(x)$ and $N_0(x)\in \mathbb{N}$ that
\begin{equation}\label{b3}
\sup_{N\geq N_0(x)} \e G(|\log \xi_{N,x}|)\leq C_0(x)<\infty.
\end{equation}

{\it Step 3} is devoted to validity of \eqref{b3}. We will employ the following statement, its proof is contained in Appendix.

\begin{Lem}\label{lemma_G}
Let $F(u)$, $u\in \mathbb{R}$, be a cumulative distribution function and $F(0)=0$. Then
\vskip0.2cm
$1)$ $\int_{(0, \frac{1}{e}]} (-\log{u}) \log{(-\log{u})} \, dF(u) = \int_{(0, \frac{1}{e}]} F(u) \frac{\log{(-\log{u})} + 1}{u} \, du$,

\vskip0.2cm
$2)$ $\int_{(e, \infty)} \log{u} \, \log{\log{u}} \, dF(u) = \int_{(e, \infty)} (1 - F(u)) \frac{\log{\log{u}} + 1}{u} \, du$.

\end{Lem}

Note that, for $u \in (\frac{1}{e},e]$, we have $|\log u|\leq 1$ and consequently
$G(|\log u|)=0$. Obviously, $F_{N,x}(0)=0$ according to \eqref{eq3a} if  $x\in \mathbb{R}^d$  and $N\geq 2$. Therefore due to Lemma \ref{lemma_G}, for such $x$ and $N$,
\begin{gather*}
    \e G(|\log{\xi_{N,x}}|) =
    \int_{(0, \frac{1}{e}]} (-\log{u}) \log{(-\log{u})} \, dF_{N,x}(u) + \int_{[e, \infty)} \log{u} \log\log{u} \, dF_{N,x}(u)  \\ = \int_{(0, \frac{1}{e}]} F_{N,x}(u) \frac{\log{(-\log{u})} + 1}{u} \, du + \int_{[e, \infty)} (1 - F_{N,x}(u)) \frac{\log{\log{u}} + 1}{u} \, du\\:=I_1(N,x)+I_2(N,x).
\end{gather*}
Consider $I_1(N,x)$. Set, for $N\geq 2$, $x\in \mathbb{R}^d$ and $u>0$,
\begin{equation}\label{p_int}
p_{N,x}(u) = \int_{B(x,r_N(u))} f(y) \, dy.
\end{equation}
Then, for  $R_1>0$ and any $u\in (0,\frac{1}{e}]$, we get
$$
r_N(u) = {\bigg( \frac{u}{\widetilde\gamma (N-1) V_d} \bigg)}^{1/d} \leq {\bigg( \frac{1}{e \widetilde\gamma (N-1) V_d} \bigg)}^{1/d} \leq R_1
$$
if $N\geq N_1$ where $N_1=N_1(R_1)$. Thus, for $R_1>0$, $u\in (0,\frac{1}{e}]$, $x\in \mathbb{R}^d$ and $N\geq N_1$,
\begin{gather*}
    \frac{p_{N,x}(u)}{|B(x, r_N(u))|} = \frac{\int_{B(x,r_N(u))} f(y) \, dy}{r_N^d(u) V_d} \leq \sup_{r \in (0, R_1]} \frac{\int_{B(x,r)} f(y) \, dy}{r^d V_d} = M_f(x,R_1)
\end{gather*}
and we obtain an inequality
\begin{equation}\label{third: eq1}
p_{N,x}(u) \leq M_f(x, R_1) \, |B(x, r_N(u))| =  \frac{M_f(x, R_1) u}{\widetilde\gamma (N-1)}.
\end{equation}
Note that if $\varepsilon \in (0,1]$ and $x \in [0,1]$ then, for all $N \geq 1$,
\begin{equation}\label{Bernoulli}
1-(1-x)^N \leq (Nx)^{\varepsilon}.
\end{equation}
Indeed, according to the Bernoulli inequality $1-(1-x)^N \leq Nx$. If $x\in [0,\frac{1}{N}]$ then
$Nx\leq 1$ and $Nx\leq (Nx)^{\varepsilon}$. If $x\in (\frac{1}{N},1]$ then $Nx >1$ and
$(Nx)^{\varepsilon}>1 \geq 1-(1-x)^N$.

By assumptions of the Theorem  $Q_f(\varepsilon_1, R_1) < \infty$ for some $\varepsilon_1 > 0$, $R_1>0$. According to Lemma \ref{lemma1} one can assume that $\varepsilon_1 < 1$. Thus, due to  \eqref{Bernoulli} and since $p_{N,x}(u) \in [0,1]$ for all $x \in \mathbb{R}^d$, $u > 0$ and $N \geq 2$, one has
\begin{equation}\label{apply_lemma3}
1-(1-p_{N,x}(u))^{N-1} \leq ((N-1)p_{N,x}(u))^{\varepsilon_1}.
\end{equation}

In view of \eqref{a1}, \eqref{p_int}, \eqref{third: eq1} and \eqref{apply_lemma3} one can claim now that, for all $x \in \Lambda(f) \cap S(f)$, $u \in (0,\frac{1}{e}]$ and $N \geq N_1$,
\vspace{-0.5cm}
\begin{align}\label{ineq_B}
\begin{gathered}
    F_{N,x}(u) = 1-(1-p_{N,x}(u))^{N-1} \leq \Big((N-1) \frac{M_f(x, R_1) u}{\widetilde\gamma (N-1)}\Big)^{\varepsilon_1} = \frac{(M_f(x, R_1))^{\varepsilon_1}}{\widetilde\gamma^{\varepsilon_1}} u^{\varepsilon_1}.
\end{gathered}
\end{align}

Therefore, for any $x\in \Lambda(f) \cap S(f)$ and $N\geq N_1$, one has
\begin{align}\label{third: eq2}
\begin{gathered}
I_1(N,x) \leq \frac{(M_f(x, R_1))^{\varepsilon_1}}{\widetilde\gamma^{\varepsilon_1}} \int_{(0, \frac{1}{e}]} \frac{\log{(-\log{u})} + 1}{u^{1-\varepsilon_1}}  \, du \\
= \frac{(M_f(x, R_1))^{\varepsilon_1}}{\widetilde\gamma^{\varepsilon_1}} \int_{[1, \infty)} (\log{t}+1) e^{-\varepsilon_1 t} dt = \frac{(M_f(x, R_1))^{\varepsilon_1}}{\widetilde\gamma^{\varepsilon_1}} L(\varepsilon_1)
\end{gathered}
\end{align}
where $L(\varepsilon) := \int_{[1, \infty)} (\log{t}+1) e^{-\varepsilon t} dt < \infty$ for each $\varepsilon > 0$.

Now consider $I_2(N,x)$. For $N\geq 10> 1+e^2$, we can write
\begin{gather*}
I_2(N,x)=
    \left(\int_{[e, {\sqrt{N-1}}]}+ \int_{({\sqrt{N-1}},{\infty})}\right) (1 - F_{N,x}(u)) \frac{\log{\log{u}} + 1}{u} \, du
    := J_1(N,x) + J_2(N,x).
\end{gather*}
For $R_2 > 0$ and any $u\in [e,\sqrt{N-1}]$, take $N_2=N_2(R_2) \geq 10$ such that if $N\geq N_2$ then
\begin{equation}\label{h1}
r_N(u) = \left( \frac{u}{\widetilde{\gamma} (N-1) V_d} \right)^{1/d} \leq \left(\frac{1}{\widetilde{\gamma}\sqrt{N-1}V_d}\right)^{1/d} \leq R_2.
\end{equation}
For the same $u$, an elementary inequality $1-t \leq e^{-t}$, $t\in [0,1]$, leads
to the following one
$$
    1 - F_{N,x}(u) = {\Big(1 -  p_{N,x}(u)  \Big)}^{N-1} \leq \exp\Big\{  -(N-1) p_{N,x}(u) \Big\}  =
        \exp\Bigg\{ -\left(\frac{u}{\widetilde\gamma}\right) \frac{p_{N,x}(u)}{\frac{u}{\widetilde\gamma (N-1)}}\Bigg\}
$$
\vspace{-0.5cm}
\begin{align}\label{1-F_add}
\begin{gathered}
        = \exp\Bigg\{ -\left(\frac{u}{\widetilde\gamma}\right) \frac{\int_{B(x,r_N(u))} f(y) \, dy}{r_N^d(u) V_d}\Bigg\}
       \leq \exp\bigg\{ -\frac{u}{\widetilde\gamma} m_f(x,R_2) \bigg\}.
\end{gathered}
\end{align}

We use an auxiliary result, its proof is provided in Appendix.
\begin{Lem}\label{l1}
For a version of a density $f$ and each $R>0$, one has $\mu(S(f)\setminus D_f(R))=0$ where $D_f(R):=\{x\in S(f): m_f(x,R)>0\}$ and $m_f(\cdot,R)$ is defined in \eqref{mM}.
\end{Lem}

It is easily seen that, for any $t>0$ and each $\delta \in (0,e]$, one has $e^{-t}\leq t^{-\delta}$.
Thus, for $x\in D_f(R_2)$, $N\geq N_2$,  $u\in [e,\sqrt{N-1}]$ and arbitrary $\varepsilon \in (0,e]$,
taking into account that $m_f(x,R_2)>0$ for $x\in D_f(R_2)$ and applying  relation \eqref{1-F_add}, we have
\begin{equation}\label{1-F}
1 - F_{N,x}(u)\leq \left(\frac{u}{\widetilde\gamma} m_f(x,R_2)\right)^{-\varepsilon}.
\end{equation}
Thus for all $x \in \Lambda(f)\cap S(f) \cap D_f(R_2)$ and any $N \geq N_2$
\begin{equation}\label{third: eq3}
    J_1(N,x)
     \leq \frac{{\widetilde\gamma}^{\varepsilon}}{{\big(m_f(x,R_2)\big)}^{\varepsilon}} \int_{[e, \infty)}
     \frac{\log{\log{u}} + 1}{u^{1+\varepsilon}} \, du
          =  \frac{{\widetilde\gamma}^{\varepsilon}}{{\big(m_f(x,R_2)\big)}^{\varepsilon}} L(\varepsilon).
\end{equation}

Further on
\begin{gather*}
    J_2(N,x)
    \leq {\Big(1 -  p_{N,x}(\sqrt{N-1})\Big)}^{N-2} \int_{({\sqrt{N-1}}, {\infty})} \frac{\log{\log{u}} + 1}{u} \Big(1 -   p_{N,x}(u)\Big) \, du \\:= J_2^{(1)}(N,x) \cdot J_2^{(2)}(N,x).
    \end{gather*}
Note that $r_N(\sqrt{N-1})\leq R_2$ for $N\geq N_2$ by virtue of \eqref{h1} and, for $N\geq 3$, one has $\frac{N-2}{N-1}\geq \frac{1}{2}$.
Hence, for  $N\geq N_2$, $\varepsilon \in (0,e]$ and
$x\in D_f(R_2)$, Lemma \ref{l1} entails
\begin{equation}\label{third: eq4}
\begin{gathered}
    J_2^{(1)}(N,x) = {\Big(1 -  p_{N,x}(\sqrt{N-1})  \Big)}^{N-2}
       \leq \exp\left\{ -(N-2) p_{N,x}(\sqrt{N-1})\right\} \\
        = \exp\Bigg\{ -\left(\frac{N-2}{N-1}\right) \left( \frac{\sqrt{N-1}}{\widetilde{\gamma}} \right) \frac{p_{N,x}(\sqrt{N-1})}{\frac{\sqrt{N-1}}{\widetilde{\gamma}(N-1)}}\Bigg\} \leq \exp\bigg\{ -\frac{\sqrt{N-1}}{2 \widetilde\gamma} m_f(x,R_2) \bigg\}
    \\ \leq \frac{(2\widetilde\gamma)^{\varepsilon}}{\big( m_f(x,R_2) \big)^{\varepsilon}{(N-1)}^{\frac{\varepsilon}{2}}}.
\end{gathered}
\end{equation}
According to \eqref{eq3a} we see that $\xi_{2,x}\stackrel{law}=\widetilde{\rho}(x,\xi):=
V_d\widetilde{\gamma}\rho^d(x,\xi)$. Consequently,
\begin{gather*}
    \e G\big( |\log{\xi_{2,x}}| \big) = \e G\big( |\log{\widetilde{\rho}(x,\xi)}| \big) = \int_{\mathbb{R}^d} {G\big(|\log{\widetilde{\rho}(x,y)}|\big)} f(y) \,dy.
\end{gather*}

\begin{Lem}\label{l4}
There are constants $a,b\geq 0$ such that, for each $x\in \mathbb{R}^d$,
\begin{equation}\label{v0}
\int\limits_{\mathbb{R}^d} {G\big(|\log{\widetilde{\rho}(x,y)}|\big)} f(y) \,dy
\leq a \int\limits_{\mathbb{R}^d} {G\big(|\log{\rho(x,y)}|\big)} f(y) \,dy + b.
\end{equation}
\end{Lem}

The proof is given in Appendix.

Set $A_f(G):=\{x\in S(f): \e G(|\log \xi_{2,x}|)<\infty\}$. By Lemma \ref{l4} one has $\mu(S(f)\setminus A_f(G))=0$ as we assumed that $K_f(\varepsilon_0)<\infty$.
Further we  write $A:=\Lambda(f)\cap S(f)\cap D_f(R_2)\cap A_f(G)$ where $\mu(S(f)\setminus A)=0$.
Introduce $w=\frac{u}{N-1}$. Then
\begin{gather*}
    J_2^{(2)}(N,x)
    = \int_{(\frac{1}{{\sqrt{N-1}}}, {\infty})} \frac{\log{\log{(w (N-1))}} + 1}{w} \Big(1 - p_{N,x}(w(N-1)) \Big)\,dw  \\
    = \bigg( \int_{(\frac{1}{{\sqrt{N-1}}}, e]} + \int_{(e, \infty)} \bigg) \frac{\log{\log{(w (N-1))}} + 1}{w} \big(1 - F_{2,x}(w) \big) \, dw  \\
    := J_{2,1}^{(2)}(N,x) + J_{2,2}^{(2)}(N,x)
\end{gather*}
where we have used that $r_N\big(w(N-1)\big) =
\Big(\frac{w}{ V_d \widetilde\gamma}\Big)^{\frac{1}{d}} = r_2(w)$
and $p_{2,x}(w)=F_{2,x}(w)$ in view of \eqref{a1}. Thus
\vspace{-0.2cm}
\begin{align}\label{third: eq5}
\begin{gathered}
J_{2,1}^{(2)}(N,x)
\leq \int_{(\frac{1}{{\sqrt{N-1}}}, e]} \big(\log{\log{(w (N-1))}} + 1 \big) \, d \log{w} \\
\leq\big(\log{\log{(e (N-1))}} + 1 \big) \big(1 + \frac{1}{2} \log{(N-1)}\big).
\end{gathered}
\end{align}
Now we  estimate $J_{2,2}^{(2)}(N)$. For $w\geq e$ and $N\geq 4$
\begin{gather*}
    \log{\log{(w(N-1))}} = \log{(\log{w} + \log{(N-1)})}  = \log{\big({\textstyle\frac{\log{w}}{\log{(N-1)}}}+1\big)} + \log{\log{(N-1)}}  \\
    \leq \log{(\log{w}+1)} + \log{\log{(N-1)}}.
\end{gather*}
Consider a function $g(w) = \frac{\log{(\log{w}+1)}}{\log{\log{w}}}$ for $w > e$. Clearly,
$$
g'(w) = \frac{1}{{w(\log{\log{w}})}^2}  \Big( \frac{\log{\log{w}}}{\log{w}+1} - \frac{\log{(\log{w}+1)}}{\log{w}} \Big) < 0, \;\;w>e.
$$
Take an arbitrary $\Delta >0$. Then, for any $w>e^{1+\Delta}$,
$$
g(w) \leq g(e^{1+\Delta}) = \frac{\log{(2+\Delta)}}{\log{(1+\Delta)}} := K(\Delta).
$$
Hence, $\log{(\log{w}+1)} \leq K(\Delta)\, \log{\log{w}}$ for $w \geq e^{1+\Delta}$ and consequently (for $N\geq 4$)
\begin{equation}\label{third: eq6}
\begin{gathered}
\log{\log{(w(N-1))}} + 1 \leq \log{(\log{w}+1)} + \log{\log{(N-1)}} + 1  \\
\leq K(\Delta) \log{\log{w}} + \log{\log{(N-1)}} + 1 \leq K(\Delta) \big( \log{\log{w}} + 1 \big) + \log{\log{(N-1)}}
\end{gathered}
\end{equation}
since $K(\Delta) > 1$. Relation \eqref{third: eq6} yields
$$
    J_{2,2}^{(2)}(N,x)
    = \bigg( \int_{({e},{e^{1+\Delta}}]} + \int_{({e^{1+\Delta}}, {\infty})} \bigg) \frac{\log{\log{(w(N-1))}} + 1}{w} \big(1 - F_{2,x}(w)  \big) \, dw
\vspace{-0.4cm}
$$
\begin{gather*}
    \leq \big( \log{\log{(e^{1+\Delta}(N-1))}} + 1 \big) \int_{({e},{e^{1+\Delta}}]} d \log{w}\\ +
    K(\Delta) \int_{({e^{1+\Delta}}, {\infty})} \frac{\log{\log{w}} + 1}{w} \big(1 - F_{2,x}(w)  \big) \, dw   \\
    + \log{\log{(N-1)}} \int_{({e^{1+\Delta}}, {\infty})} \frac{1}{w} \big(1 - F_{2,x}(w)  \big) \, dw
    \leq \Big( \log{\log{(e^{1+\Delta} (N-1))}} + 1 \Big) \Delta  \\
    + K(\Delta) \int_{(e, \infty)} \frac{\log{\log{w}} + 1}{w} \big(1 - F_{2,x}(w)  \big) \, dw +
    \log{\log{(N-1)}} \int_{(e, \infty)} \frac{1}{w} \big(1 - F_{2,x}(w)  \big) \, dw.
\end{gather*}
By Lemma \ref{lemma_G}
$$
\int_{(e, \infty)} \frac{1}{w} \big(1 - F_{2,x}(w)  \big) \, dw \leq \int_{(e, \infty)} \frac{\log{\log{w}} + 1}{w} \big(1 - F_{2,x}(w)  \big) \, dw
$$
$$
= \int_{(e, \infty)} \log(w) \log{\log(w)} dF_{2,x}(w) \leq \e G(|\log \xi_{2,x}|):=R(x)<\infty
$$
whenever $x\in A$. Therefore
\begin{equation}\label{third: eq7}
\begin{gathered}
    J_{2,2}^{(2)}(N,x) \leq \Big( \log{\log{(e^{1+\Delta} (N-1))}} + 1 \Big)\Delta + K(\Delta) R(x) + \log{\log{(N-1)}} R(x).
\end{gathered}
\end{equation}
For $x\in A$ and $N\geq N_2$, relations \eqref{third: eq4}, \eqref{third: eq5}, \eqref{third: eq7} imply
\begin{equation}\label{third: eq8}
\begin{gathered}
    J_2(N,x) = J_2^{(1)}(N,x)  (J_{2,1}^{(2)}(N,x) + J_{2,2}^{(2)}(N,x))\\
    \leq \frac{(2\widetilde\gamma)^{\varepsilon}}{{(N-1)}^{\frac{\varepsilon}{2}}\big( m_f(x,R_2) \big)^{\varepsilon} }  \Big(\big(\log{\log{(e (N-1))}} + 1 \big) \big(1 + \frac{1}{2} \log{(N-1)}\big)  \\ + \big( \log{\log{(e^{1+\Delta} (N-1))}} + 1 \big) \Delta + K(\Delta) R(x) + \log{\log{(N-1)}} R(x) \Big):= R_f(N,x,\Delta).
\end{gathered}
\end{equation}
Thus, for each $\kappa >0$, $x\in A$ and $N>N_3(x)$, we can claim that $J_2(N,x)<\kappa$. In such a way, for any $N \geq \max\{N_2,N_3(x)\}:=N_4(x)$, taking into account \eqref{third: eq3} one infers that
\begin{gather*}
    I_2(N,x)
    = J_1(N,x) + J_2(N,x)  \leq \frac{L(\varepsilon) {\widetilde\gamma}^{\varepsilon}}{{(m_f(x,R_2))}^{\varepsilon}} + \kappa < \infty.
\end{gather*}

Due to Lemma \ref{lemma1} we can assume that $\varepsilon_2 \leq e\frac{1+\varepsilon_0}{\varepsilon_0}$ in conditions of Theorem 1.
Choose  $\varepsilon := \frac{\varepsilon_0\varepsilon_2}{1+\varepsilon_0}$ then
$\varepsilon \leq \frac{\varepsilon_0}{1+\varepsilon_0}e \frac{1+\varepsilon_0}{\varepsilon_0}\leq e$.
Thus, for $x\in A$ and $N\geq \max\{N_1,N_4(x)\} := N_0(x)$ and $\varepsilon = \frac{\varepsilon_0\varepsilon_2}{1+\varepsilon_0}$, according to \eqref{third: eq2}
\begin{gather*}
\e G(|\log{\xi_{N,x}}|) = I_1(N,x) + I_2(N, x) \\ \leq L(\varepsilon_1)\widetilde\gamma^{-\varepsilon_1} (M_f(x, R_1))^{\varepsilon_1}  + L(\varepsilon) {\widetilde\gamma}^{\varepsilon} {(m_f(x,R_2))}^{-\varepsilon} + \kappa := C_0(x)
\end{gather*}
and we come to the desired relation \eqref{b3}. Hence relation \eqref{convc} holds as well.

\vskip0.3cm
{\it Step 4}. Now we are able to verify \eqref{main1}. We have already proved for each $x\in A$ (thus, for $\p_{\xi}$-almost every $x$ into the support $S(f)$) that $\e(\zeta_1(N)|X_1=x) \to -\log f(x), N\to \infty$.

Set $Y_N(x):=\e(\zeta_1(N)|X_1=x)$. We have seen that $\e(\zeta_1(N)|X_1=x)=\e \log \xi_{N,x}$.
Consider $x\in A$ and take any $N\geq \max\{N_1,N_2\}$. Function $G$ is nondecreasing and convex.
On account of the Jensen inequality
\begin{equation}\label{fourth: eq1}
\begin{gathered}
    G(|Y_N(x)|) = G(|\e\log{\xi_{N,x}}|) \leq G(\e |\log{\xi_{N,x}}|)
     \leq  \e G(|\log{\xi_{N,x}}|).
\end{gathered}
\end{equation}
Relations \eqref{third: eq2}, \eqref{third: eq3}, \eqref{third: eq8} and \eqref{fourth: eq1} guarantee that, for each $x\in A$ and all $N \geq \max\{N_1,N_2\}$,
\begin{gather*}
    G(|Y_N(x)|) \leq \frac{L(\varepsilon_1) (M_f(x, R_1))^{\varepsilon_1}}{\widetilde\gamma^{\varepsilon_1}}  + \frac{L(\varepsilon) {\widetilde\gamma}^{\varepsilon}}{\big( m_f(x,R_2) \big)^{\varepsilon}}
    + R_f(N,x,\Delta).
    \end{gather*}

Thus, for any $N\geq \max\{N_1,N_2\}$,
\begin{align}\label{fourth: eq2}
\begin{gathered}
    \int_{\mathbb{R}^d} G(|Y_N(x)|) f(x) \, dx = \int_{A} G(|Y_N(x)|) f(x) \, dx  \\ \leq \frac{L(\varepsilon_1)}{\widetilde\gamma^{\varepsilon_1}} \int_{A} M_f^{\varepsilon_1}(x, R_1) f(x) \, dx + L(\varepsilon) {\widetilde\gamma}^{\varepsilon} \int_{A} m_f^{-\varepsilon}(x,R_2) f(x) \, dx  \\ + \frac{(2\widetilde\gamma)^{\varepsilon}}{{(N-1)}^{\frac{\varepsilon}{2}}}  \bigg\{ \big(\log{\log{(e (N-1))}} + 1 \big) \big(1 + \frac{1}{2} \log{(N-1)}\big)\\ + \big( \log{\log{(e^{1+\Delta} (N-1))}} + 1 \big) \Delta \bigg\}  \int_{A} m_f^{-\varepsilon}(x,R_2) f(x) \, dx  \\ + \frac{(2\widetilde\gamma)^{\varepsilon}}{{(N-1)}^{\frac{\varepsilon}{2}}} \bigg\{ K(\Delta) + \log{\log{(N-1)}}\bigg\}  \int_{A} R(x) \, m_f^{-\varepsilon}(x,R_2) f(x) \, dx.
\end{gathered}
\end{align}
Clearly,
\begin{gather}\label{u0}
\int_{A} M_f^{\varepsilon_1}(x, R_1) f(x) \, dx = \int_{\mathbb{R}^d} M_f^{\varepsilon_1}(x, R_1) f(x) \, dx = Q_f(\varepsilon_1, R_1) < \infty
\end{gather}
according to Theorem \ref{th1} assumptions.

Recall that $\varepsilon = \frac{\varepsilon_0\varepsilon_2}{1+\varepsilon_0}$. By the Lyapunov inequality one has
\begin{align}\label{u1}
\begin{gathered}
\int_{A} m_f^{-\varepsilon}(x,R_2) f(x) \, dx
\leq {\bigg( \int_{\mathbb{R}^d} m_f^{-\varepsilon \frac{1+\varepsilon_0}{\varepsilon_0}}(x,R_2) f(x) \, dx \bigg)}^{\frac{\varepsilon_0}{1+\varepsilon_0}} = \big(T_f(\varepsilon_2,R_2)\big)^{\frac{\varepsilon_0}{1+\varepsilon_0}} < \infty.
\end{gathered}
\end{align}
Now the H\"older inequality, relation \eqref{v0} and the conditions $T_f(\varepsilon_2,R_2)<\infty$, $K_f(\varepsilon_0)<\infty$ yield
\vspace{-0.3cm}
\begin{align}\label{u2}
\begin{gathered}
\int_{A} R(x) \, m_f^{-\varepsilon}(x,R_2) f(x) \, dx \\
\leq
\bigg(\int_{\mathbb{R}^d} R^{1+\varepsilon_0}(x) f(x) \, dx \bigg)^{\frac{1}{1+\varepsilon_0}}
{\bigg( \int_{\mathbb{R}^d} m_f^{-\varepsilon \frac{1+\varepsilon_0}
{\varepsilon_0}}(x,R_2) f(x) \, dx \bigg)}^{\frac{\varepsilon_0}{1+\varepsilon_0}} \\
= \Big(\int\limits_{\mathbb{R}^d} {\Big(\int\limits_{\mathbb{R}^d} {G\big(|\log{\widetilde{\rho}}(x,y)|\big)} f(y) \,dy \Big)}^{1+\varepsilon_0}\!\! f(x) \, dx\Big)^{\frac{1}{1+\varepsilon_0}} \big(T_f(\varepsilon_2,R_2)\big)^{\frac{\varepsilon_0}{1+\varepsilon_0}} \\
\leq
\Big(\int\limits_{\mathbb{R}^d}\!\! {\Big(a \int\limits_{\mathbb{R}^d} {G\big(|\log{\rho}(x,y)|\big)} f(y) \,dy + b\Big)}^{1+\varepsilon_0}\!\! f(x) \, dx\Big)^{\frac{1}{1+\varepsilon_0}} \big(T_f(\varepsilon_2,R_2)\big)^{\frac{\varepsilon_0}{1+\varepsilon_0}} \\
\leq
2^{\frac{\varepsilon_0}{1+\varepsilon_0}}(a \big(K_f(\varepsilon_0)\big)^{\frac{1}{1+\varepsilon_0}}+b)\, \big(T_f(\varepsilon_2,R_2)\big)^{\frac{\varepsilon_0}{1+\varepsilon_0}}
< \infty.
\end{gathered}
\end{align}
where $a, b \geq 0$ are the same as in Lemma \ref{l4}.
Accordingly, for $N\geq \max\{N_1,N_2\}$ we can rewrite \eqref{fourth: eq2} in the following way
\begin{gather*}
    \int_{\mathbb{R}^d} G(|Y_N(x)|) f(x) \, dx \leq
         S_1 + S_2 + S_3(N) + S_4(N)
\end{gather*}
where $S_3(N) \rightarrow 0$, $S_4(N) \rightarrow 0$ as $N \rightarrow \infty$. Thus,
\begin{gather*}
   \sup_{N\geq \max\{N_1,N_2\}}\int_{\mathbb{R}^d} G(|Y_N(x)|) f(x) \, dx  < \infty.
\end{gather*}
We have established  uniform integrability of  the family $\{Y_N\}_{N \geq \max\{N_1, N_2\}}$
(w.r.t. measure $\p_{\xi}$). Therefore, we conclude that
\begin{gather*}
    \e\zeta_1(N) = \int_{\mathbb{R}^d} Y_N(x) f(x) \, dx \to -\int_{\mathbb{R}^d}  \log{f(x)}f(x) \, dx = H,\;\;N\to \infty.
\end{gather*}
So we come to \eqref{main1}. The proof of Theorem 1 is complete. $\square$
\vskip0.2cm

{\it Proof of Corollary 2}.
Assume that condition (A) is satisfied. In Section 2 we used the elementary inequality:
$\log z \leq \frac{1}{\delta}z^{\delta}$ for any $z\geq 1$ and $\delta >0$.
Hence, for any $t \geq 0$ and $\delta >0$, $G(t)\leq \frac{1}{\delta}t^{1+\delta}$.
By the Lyapunov inequality, for any $\delta>0$ and $\varepsilon_0 >0$,
$$
\left(\int_{\mathbb{R}^d}|\log\rho(x,y)|^{1+\delta}
f(y)\,dy\right)^{1+\varepsilon_0}\leq \int_{\mathbb{R}^d}|\log\rho(x,y)|^{(1+\delta)(1+\varepsilon_0)}
f(y)\,dy.
$$
For a given $p>1$ take $\delta:=\frac{p-1}{p+1}$ and $\varepsilon_0:=\frac{p-1}{2}$ then
$(1+\delta)(1+\varepsilon_0)=p$. Thus
\begin{gather*}
K_f(\varepsilon_0)  \leq
\left(\frac{1}{\delta}\right)^{1+\varepsilon_0} \int_{\mathbb{R}^d}\int_{\mathbb{R}^d}|\log\rho(x,y)|^{p}
f(x) f(y)\,dx\,dy < \infty.
\end{gather*}
Hence condition (A) implies that $K_f(\varepsilon_0)<\infty$ for $\varepsilon_0=\frac{p-1}{2}$.

If (B) is true then, for each $x\in \mathbb{R}^d$ and any $R>0$, one has $M_f(x,R)\leq M$.
Therefore $Q_f(\varepsilon_1, R_1)<\infty$ for any $\varepsilon_1>0$ and $R_1>0$.

Assume that (C1) is satisfied. Accordingly $m_f(x,R)\geq m>0$ for any $x\in \mathbb{R}^d$ and $R>0$.
Then evidently $T_f(\varepsilon_2,R_2)<\infty$ for any $\varepsilon_2>0$ and $R_2>0$.

It remains to show that if (B) and (C1) are satisfied then \eqref{main1} holds whenever $f$ has a bounded support. It is sufficient to verify that if $f\leq M$ and $f$ has a bounded support then (A) is valid. Clearly, for any ball $B(0,r)\subset \mathbb{R}^d$ such that $S(f)\subset B(0,r)$,
$$
\int_{\mathbb{R}^d}\int_{\mathbb{R}^d}|\log\rho(x,y)|^{p}
f(x) f(y)\,dx\,dy \leq M^2\int_{B(0,r)}\int_{B(0,r)}|\log\rho(x,y)|^{p}\,dx\, dy.
$$
Note that $\mu(B(0,r))<\infty$. Changing the variables of integration in inner integral over $x$ (now we use the spherical coordinates) and taking into account that, for any $p>1$,
$$
\int_{(0,2r)}|\log u|^p u^{d-1}du <\infty
$$
we see that (A) is true. $\square$

\section{Proof of Theorem 2}

First of all note that, for any $t>0$, $G(|\log t|)\leq G(\log^2 t)$. Indeed,
if $t\in (0,\frac{1}{e}] \cup (e,\infty)$ then $|\log t|\geq 1$. Thus, for such $t$,
$|\log t|\leq \log^2 t$ and the desired inequality holds since $G$ is a nondecreasing function.
If $t\in (\frac{1}{e},e]$ then $G(|\log t|)=G(\log^2 t)=0$.
Consequently if $K_{f,2}(\varepsilon_0)<\infty$ then $K_f(\varepsilon_0)<\infty$ and
the conditions of Theorem \ref{th_main2} guarantee the statement of Theorem \ref{th1} validity.

We will show that Theorem 2 conditions permit to claim that $\zeta_i(N) \in L^2$ for $i=1,\ldots,N$ and $N \geq N_0$ where $N_0 \in \mathbb{N}$. Then we can write $\e (H_N - H)^2 = \var(H_N) + (\e H_N-H)^2$ and
accordingly we come to \eqref{main2} if $\var (H_N)\to 0$, $N\to \infty$. Taking into account that $\zeta_1(N),\ldots,\zeta_N(N)$ are identically distributed one has
\begin{equation}\label{t2: second: eq2}
\var(H_N) =  \frac{1}{N}\var(\zeta_1(N)) + \frac{2}{N^2}\sum\limits_{1\leq i< j\leq N} \cov(\zeta_i(N), \zeta_j(N)).
\end{equation}
We do not strictly adhere to the notation used in
Theorem~\ref{th1} proof. Namely, the choice of positive constants $C_k$, $C_k(x)$ and
integers $N_k$, $N_k(x)$ where
$k\in \{0\}\cup \mathbb{N}$ and $x\in \mathbb{R}^d$ could be different.
The proof of Theorem \ref{th_main2} is also divided into several steps.

\textit{Step 1.} We study $\e \zeta_1(N)^2$. Let us formulate the following result.
\begin{Lem}\label{lemma6}
There are constants $a,b\geq 0$ such that, for each $x\in \mathbb{R}^d$,
\begin{equation*}
\int\limits_{\mathbb{R}^d} {G\big(\log^2{\widetilde{\rho}(x,y)}\big)} f(y) \,dy
\leq a \int\limits_{\mathbb{R}^d} {G\big(\log^2{\rho(x,y)}\big)} f(y) \,dy + b.
\end{equation*}
\end{Lem}

The proof is quite similar to that of Lemma \ref{l4} and therefore is skipped.

Let $A:=\Lambda(f)\cap S(f)\cap D_f(R_2)\cap A_{f,2}(G)$ where the first three sets appeared in the Theorem \ref{th1} proof and $A_{f,2}(G):=\{x\in S(f): \e G(\log^2 \xi_{2,x})<\infty\}$. It is easily seen that $\p_{\xi}(A) = 1$ because
$\p_{\xi}(\Lambda(f)\cap S(f)\cap D_f(R_2))=1$ and $\p_{\xi}(A_{f,2}(G))=1$ (we take into account the result of Lemma \ref{lemma6} and the fact that $K_{f,2}(\varepsilon_0)<\infty$).
In view of \eqref{eq1a},  for each $x \in A$,
\begin{gather*}
    \e (\zeta_1(N)^2 | X_1 = x) = \e \big( (\log e^{\zeta_1(N)})^2 | X_1 = x\big) = \int_{(0, \infty)} \log^2(u) \, dF_{N,x}(u) = \e \log^2 \xi_{N,x}.
\end{gather*}
Since $\xi_x\sim Exp\big(\frac{f(x)}{\widetilde{\gamma}}\big)$, $x\in S(f)$, one has
\begin{gather*}
    \e \log^2 \xi_x  = \int_{(0, \infty)} \log^2(u) \, d F_x(u)
    = \log^2{f(x)} + \sigma^2
\end{gather*}
where
$$
\sigma^2 := \int_{(0, \infty)} \log^2{v} \, e^{-v} \, dv - \left(\int_{(0, \infty)} \log{v} \, e^{-v} \, dv\right)^2<\infty.
$$

Due to \eqref{claw} we have
$(\log \xi_{N,x})^2 \stackrel{law}\rightarrow (\log \xi_x)^2$ for $x\in A$ as $N\to \infty$.
Now we show that for $x\in A$ a family $\left\{\log^2 \xi_{N,x}\right\}_{N \geq N_0(x)}$ is uniformly integrable for some $N_0(x) \in \mathbb{N}$. It is sufficient to verify that, for all $x\in A$,
\begin{equation}\label{t2: second: eq3}
\sup\limits_{N \geq N_0(x)} \e G(\log^2 \xi_{N,x}) \leq C_0(x) < \infty.
\end{equation}
To show \eqref{t2: second: eq3} we employ the following result.

\begin{Lem}\label{le7}
Let $F$ be a cumulative distribution function such that $F(0) = 0$. Then the following relations hold:

\vskip0.2cm

1) $\int_{(0, \frac{1}{e}]} \log^2{u} \log(-\log{u}) \, dF(u) = 2 \int_{(0, \frac{1}{e}]} F(u) \frac{(-\log{u})}{u} \big( \log{(-\log{u})} + \frac{1}{2} \big) \, du$,

\vskip0.1cm

2) $\int_{[e, \infty)} \log^2{u} \log{\log{u}} \, dF(u) = 2 \int_{[e, \infty)} (1 - F(u)) \frac{\log{u}}{u} \big( \log{\log{u}} + \frac{1}{2} \big) \, du$.

\end{Lem}

The proof is quite similar to that of Lemma \ref{lemma_G} and therefore is skipped.
By Lemma \ref{le7}
\begin{align}\label{t2: second: eq4}
\begin{gathered}
    \e G(\log^2{\xi_{N,x}})
     =
     \int_{(0, \frac{1}{e}]} G(\log^2{u}) \, dF_{N,x}(u) + \int_{(e, \infty)} G(\log^2{u}) \, dF_{N,x}(u)  \\ =
    2 \int_{(0, \frac{1}{e}]} \log^2{u} \log{(-\log{u})} \, dF_{N,x}(u) + 2 \int_{(e, \infty)} \log^2{u} \log\log{u} \, dF_{N,x}(u)  \\ = 4 \int_{(0, \frac{1}{e}]} F_{N,x}(u) \frac{(-\log{u})}{u} \big( \log{(-\log{u})} + \frac{1}{2} \big) \, du \\ + 4 \int_{(e, \infty)} (1 - F_{N,x}(u)) \frac{\log{u}}{u} \big( \log{\log{u}} + \frac{1}{2} \big) \, du\\:=4 I_1(N,x)+4 I_2(N,x).
\end{gathered}
\end{align}

Employment of \eqref{ineq_B} leads, for $x\in A$ and $N \geq N_1$ ($N_1$ is the same as in the proof of Theorem \ref{th1}),
to relations
\begin{gather*}
I_1(N,x) = \int_{(0, \frac{1}{e}]} F_{N,x}(u) \frac{(-\log{u})}{u} \big( \log{(-\log{u})} + \frac{1}{2} \big) \, du  \\ \leq \frac{(M_f(x, R_1))^{\varepsilon_1}}{\widetilde\gamma^{\varepsilon_1}} \int_{(0, \frac{1}{e}]} \frac{(-\log{u})
( \log{(-\log{u})} + \frac{1}{2} )}{u^{1-\varepsilon_1}} \, du = M_f^{\varepsilon_1}(x, R_1) L_1(\varepsilon_1, \widetilde\gamma)
\end{gather*}
where $L_1(\varepsilon_1, \widetilde\gamma) = \frac{1}{\widetilde\gamma^{\varepsilon_1}} \int_{[1, \infty)} v
(\log{v} + \frac{1}{2}) e^{-\varepsilon_1  v}\, dv < \infty$. Furthermore, one has
\begin{gather*}
I_2(N,x) = \left( \int_{(e, \sqrt{N-1}]}+ \int_{(\sqrt{N-1},\infty)} \right) \big(1-F_{N,x}(u)\big) \frac{\log{u}}{u} \big( \log{\log{u}} + \frac{1}{2} \big) \, du \\ := J_1(N,x) + J_2(N,x).
\end{gather*}

Set $\varepsilon := \frac{\varepsilon_0\varepsilon_2}{1+\varepsilon_0}$ (as in the proof of Theorem \ref{th1} w.l.g. we can assume that $\varepsilon \leq e$). For $N \geq N_2$ ($N_2$ is the same as in the proof of Theorem \ref{th1}), taking into account inequality \eqref{1-F} we get
\begin{gather*}
J_1(N,x) \leq \frac{{\widetilde\gamma}^{\varepsilon}}{{\big(m_f(x,R_2)\big)}^{\varepsilon}} \int_{(e, \infty)} \frac{\log{u}(\log{\log{u}} + \frac{1}{2})}{u^{1+\varepsilon}} \, du = m_f^{-\varepsilon}(x,R_2) L_2(\varepsilon, \widetilde\gamma)
\end{gather*}
where $L_2(\varepsilon, \widetilde\gamma) = \widetilde\gamma^{\varepsilon} \int_{[1,\infty)} v (\log{v} + \frac{1}{2}) e^{-\varepsilon v} \, dv < \infty$.
Since $(1-F_{N,x}(u)) = \big( 1 - p_{N,x}(u) \big)^{N-1}$ we infer that
\begin{gather*}
    J_2(N,x) \leq \big( 1 - p_{N,x}(\sqrt{N-1})\big)^{N-2} \int_{(\sqrt{N-1}, \infty)} \frac{\log{u} (\log{\log{u}} + \frac{1}{2})}{u} (1- p_{N,x}(u))\, du.
\end{gather*}
\vspace{-0.5cm}

By employing \eqref{third: eq4} and \eqref{third: eq6} we deduce, for $w = \frac{u}{N-1}$, the inequality
$$
    J_2(N,x) \leq \frac{(2\widetilde\gamma)^{\varepsilon}}{\big( m_f(x,R_2) \big)^{\varepsilon}{(N-1)}^{\frac{\varepsilon}{2}}} \times
$$
$$
    \times \left( \int_{(\frac{1}{\sqrt{N-1}}, e^{1+\Delta}]} + \int_{(e^{1+\Delta}, \infty)} \right) \frac{\log{(w (N-1))} (\log{\log{(w (N-1))}} + \frac{1}{2})}{w} \big(1 - F_{2,x}(w) \big) \, dw
$$
\vspace{-0.3cm}
$$
    \leq \frac{(2\widetilde\gamma)^{\varepsilon}}{\big( m_f(x,R_2) \big)^{\varepsilon}{(N-1)}^{\frac{\varepsilon}{2}}} \times
$$
$$
    \times \bigg\{ \log(e^{1+\Delta} (N-1))\big(\log{\log{(e^{1+\Delta} (N-1))}} + \frac{1}{2} \big) \big(1 + \Delta + \frac{1}{2} \log{(N-1)}\big) + \log{(N-1)} \times \\
$$
$$
    \times \Big[ K(\Delta) \int_{(e, \infty)} \frac{\log{\log{w}} + \frac{1}{2}}{w} \big(1 - F_{2,x}(w)  \big) \, dw +
    \log{\log{(N-1)}} \int_{(e, \infty)} \frac{1}{w} \big(1 - F_{2,x}(w)  \big) \, dw \Big]
$$
\vspace{-0.5cm}
\begin{align}\label{J_2_expand}
\begin{gathered}
    + K(\Delta) \int_{(e, \infty)} \frac{\log{w} (\log{\log{w}} + \frac{1}{2})}{w} \big(1 - F_{2,x}(w)  \big) \, dw \\ +
    \log{\log{(N-1)}} \int_{(e, \infty)} \frac{\log{w}}{w} \big(1 - F_{2,x}(w)  \big) \, dw  \bigg\}
\end{gathered}
\end{align}
where $\Delta > 0$ is an arbitrary number and $K(\Delta) = \frac{\log{(2+\Delta)}}{\log{(1+\Delta)}}$.

Set $P(x) := \e G(\log^2 \xi_{2,x})$, $x \in A_{f,2}(G)$. Then $\int_{[e, \infty)} (1-F_{2,x}(w)) \frac{\log{w} (\log\log{w} + \frac{1}{2})}{w} \, dw \leq P(x) < \infty$
in view of Lemma \ref{le7}. Each integrals appearing in \eqref{J_2_expand} admits the same bound (up to a constant factor). Therefore, for all $x\in A$ and $N \geq N_2$
\begin{gather*}
    J_2(N,x) \leq m_f^{-\varepsilon}(x,R_2) \big(h_1(N,\varepsilon, \widetilde\gamma) + P(x) h_2(N,\varepsilon, \widetilde\gamma)\big).
\end{gather*}
where $h_1(N,\varepsilon, \widetilde\gamma) \to 0$ and $h_2(N,\varepsilon, \widetilde\gamma) \to 0$ as $N \to \infty$.
Thus, according to \eqref{t2: second: eq4}, for any $x \in A$, $N \geq \max\{N_1, N_2\}$ and $\varepsilon = \frac{\varepsilon_0 \varepsilon_2}{1+ \varepsilon_0}$, we see that
\begin{align}\label{idea}
\begin{gathered}
    \e G(\log^2 \xi_{N,x})  \leq C_1 M_f^{\varepsilon_1}(x, R_1) + m_f^{-\varepsilon}(x, R_2)(C_2  + C_3(N) + C_4(N)  P(x))
\end{gathered}
\end{align}
where $C_1$, $C_2$, $C_3(N)$ and $C_4(N)$ can  depend only on $\varepsilon_0$, $\varepsilon_1$, $\varepsilon_2$ and $\widetilde\gamma$ ($\Delta$ is fixed). Moreover, $C_3(N) \to 0$, $C_4(N) \to 0$ as $N \to \infty$.
Hence, for any $\kappa > 0$ and $N \geq N_0(x):=\max\left\{ N_1, N_2, N_3(x) \right\}$,
\begin{gather*}
\sup_{N \geq N_0(x)} \e G(\log^2 \xi_{N,x}) \leq C_1 M_f^{\varepsilon_1}(x, R_1) + C_2 m_f^{-\varepsilon}(x, R_2) + \kappa := C_0(x).
\end{gather*}
Therefore, a family $\left\{\log^2 \xi_{N,x}\right\}_{N \geq N_0(x)}$ is uniformly integrable. Hence for any $x\in A$,
$$
\e(\zeta_1(N)^2|X_1=x)\to \log^2f(x) + \sigma^2,\;\;N\to \infty.
$$
Set now $Z_N(x):=\e(\zeta_1(N)^2|X_1=x) \geq 0$. We have already seen that $\e \log^2 \xi_{N,x}=\e(\zeta_1(N)^2|X_1=x)$. Also it was proved that $Z_N(x)\to \log^2f(x) + \sigma^2$
for each $x\in A$ as $N\to \infty$ (thus, for $\p_{\xi}$-almost every $x$ into the support $S(f)$).
By the Jensen inequality
$G(Z_N(x)) = G(\e \log^2{\xi_{N,x}}) \leq \e  G(\log^2{\xi_{N,x}})$
because $G(\cdot)$ is convex.
Then according to \eqref{idea}, for any $x\in A$ and for any $N\geq \max\{N_1,N_2\}$, one has
\begin{gather*}
G(Z_N(x)) \leq C_1 M_f^{\varepsilon_1}(x, R_1) + C_2 m_f^{-\varepsilon}(x, R_2) + C_3(N) m_f^{-\varepsilon}(x,R_2) + C_4(N) m_f^{-\varepsilon}(x, R_2) P(x).
\end{gather*}
Thus, for all $N\geq \max\{N_1,N_2\}$, in view of \eqref{u0}, \eqref{u1} and analogue of \eqref{u2} for $P(x)$ instead of $R(x)$
\begin{gather*}
  \int_{\mathbb{R}^d} G(Z_N(x)) f(x) \, dx \leq C_1 \int_{\mathbb{R}^d} M_f^{\varepsilon_1}(x,R_1) f(x) \, dx + C_2 \int_{\mathbb{R}^d} m_f^{-\varepsilon}(x, R_2) f(x) \, dx \\ + C_3(N)
  \int_{\mathbb{R}^d} m_f^{-\varepsilon}(x,R_2) f(x) \, dx + C_4(N) \int_{\mathbb{R}^d} m_f^{-\varepsilon}(x,R_2) P(x) f(x) \, dx \\ \leq \widetilde{S}_1 + \widetilde{S}_2 + \widetilde{S}_3(N) + \widetilde{S}_4(N),
\end{gather*}
where $\widetilde{S}_3(N) \to 0$, $\widetilde{S}_4(N) \to 0$ as $N \to \infty$. Thus
\vspace{-0.3cm}
$$
\sup_{N \geq \max\{N_1, N_2\}} \int_{\mathbb{R}^d} G(Z_N(x)) f(x) \, dx < \infty.
\vspace{-0.3cm}
$$
Consequently, a family $\{Z_N(\cdot)\}_{N\geq \max\{N_1,N_2\}}$ is uniformly integrable and we can claim that
\begin{equation*}
\e \zeta_1^2(N) \to \int_{\mathbb{R}^d} f(x) \log^2{f(x)} \, dx +  \sigma^2.
\end{equation*}
Hence
$\var(\zeta_1(N)) = \e\zeta_1^2(N) - (\e\zeta_1(N))^2 \to
\sigma^2 + \int_{\mathbb{R}^d} f(x) \log^2{f(x)} \, dx - H^2$
as $N\to \infty$.
As in Remark \ref{rem2} one can prove in a similar way that
finiteness of integrals in \eqref{p2} and \eqref{p3} implies that $\int_{\mathbb{R}^d} f(x) \log^2{f(x)} \, dx < \infty$.
Clearly, we have inferred that $\frac{1}{N}\var(\zeta_1(N))\to 0$ as $N\to \infty$.

\vskip0.2cm
\textit{Step 2.} Now we consider $\cov(\zeta_i(N),\zeta_j(N))$ for $i\neq j$, where $i,j\in \{1,\ldots,N\}$.
For $x\in A$ and $y\in A$ (where $A$ has been defined at the beginning of \textit{Step 1}), introduce  conditional cumulative distribution function
\begin{equation*}
    \widetilde{F}^{i,j}_{N,x,y}(u,w) := P(e^{\zeta_i(N)} \leq u, e^{\zeta_j(N)} \leq w | X_i = x, X_j = y), \;\;u,w\geq 0.
\end{equation*}
Here $i,j\in \{1,\ldots,N\}$, $i\neq j$.
For any events $C$ and $D$, one can write
$
\ind\{C\cap D\}= 1- \ind\{\overline{C}\}-\ind\{\overline{D}\}+ \ind\{\overline{C}\cap \overline{D}\}
$
where $\overline{C}:=\Omega\setminus C$. Thus, for $x,y\in \mathbb{R}^d$, $u,w \geq 0$, $i\neq j$ and $N\geq 3$,
\begin{align}\label{t2: third: eq1}
\begin{gathered}
    \widetilde{F}^{i,j}_{N,x,y}(u,w) = 1 - P(e^{\zeta_i(N)} > u | X_i = x, X_j = y) - P(e^{\zeta_j(N)} > w | X_i = x, X_j = y)  \\ + P(e^{\zeta_i(N)} > u, e^{\zeta_j(N)} > w | X_i = x, X_j = y) =
    \\ 1- P(\min_{k \neq i} \rho(X_i, X_k) > r_N(u) | X_i = x, X_j = y) \\ - P(\min_{k \neq j} \rho(X_j, X_k) > r_N(w) | X_i = x, X_j = y)  \\ + P(\min_{k \neq i} \rho(X_i, X_k) > r_N(u), \, \min_{k \neq j} \rho(X_j, X_k) > r_N(w) | X_i = x, X_j = y)\\
    = 1 - \ind{\big(\rho(x,y) > r_N(u)\big)}\,P\big(\min_{k \neq i, j} \rho(x, X_k) > r_N(u) \big)  \\ -
    \ind{\big(\rho(x,y) > r_N(w) \big)}\,P\big(\min_{k \neq i, j} \rho(y, X_k) > r_N(w)\big) \\ + \ind{\big(\rho(x,y) > \max\{r_N(u), r_N(w)\}\big)}\, P\big(\min_{k \neq i, j} \rho(x, X_k) > r_N(u), \min_{k \neq i, j} \rho(y, X_k) > r_N(w) \big)
\end{gathered}
\end{align}
where $k\in \{1,\ldots,N\}$ and, as before in \eqref{rN}, $r_N(a):=\left(\frac{a}{V_d\widetilde{\gamma}(N-1)}\right)^{\frac{1}{d}}$, $a\geq 0$.
One can write $\widetilde{F}_{N,x,y}(u,w)$ because $\widetilde{F}^{i,j}_{N,x,y}(u,w)$ is the same for all $i,j\in \{1,\ldots,N\}$, $i\neq j$.

Set $A_1 := \big\{(x,y): x \in A, \, y \in A, \, x \neq y\big\}$ and $A_2 := \big\{(x,y): x \in A, \, y \in A, \, x = y\big\}$. Evidently, $\left(\p_{\xi}\otimes \p_{\xi}\right)(A_1) = 1$ and $\left(\p_{\xi}\otimes \p_{\xi}\right)(A_2) = 0$.
In case of $(x,y) \in A_2$, one has
$$
\widetilde{F}_{N,x,y}(u,w) \equiv 1, \;\; u \geq 0, \, w\geq 0.
$$

Further we consider $(x,y) \in A_1$. For such $(x,y)$ and any $a>0$,
$r_N(a)\to 0$  as $N\to \infty$. Hence
we can find $N_0=N_0(x,y,u,w)$ such that $r_N(u)<\rho(x,y)$, $r_N(w)<\rho(x,y)$ and $B(x,r_N(u)) \cap B(y, r_N(w)) = \varnothing$ for $N\geq N_0$. Then, according to \eqref{t2: third: eq1} for these $N \geq N_0$ one has
$$
    \widetilde{F}_{N,x,y}(u,w) = 1 - P\big(\min_{k \neq i, j} \rho(x, X_k) > r_N(u)\big) -
    P\big(\min_{k \neq i, j} \rho(y, X_k) > r_N(w) \big)
$$
$$
    + P\big(\min_{k \neq i, j} \rho(x, X_k) > r_N(u), \min_{k \neq i, j} \rho(y, X_k) > r_N(w) \big)
    $$
\begin{gather*}
= 1 - \prod_{k \neq i, j} P\big( X_k \notin B(x, r_N(u)) \big) - \prod_{k \neq i, j} P\big( X_k \notin B(y, r_N(w)) \big) \\ + \prod_{k \neq i, j} P\big( X_k \notin B(x, r_N(u)) \sqcup B(y, r_N(w)) \big) = 1 - \Big( 1 - \int\limits_{B(x,r_N(u))} f(z) \, dz \Big)^{N-2}  \\ - \Big( 1 - \int\limits_{B(y,r_N(w))} f(z) \, dz \Big)^{N-2} + \Big( 1 - \int\limits_{B(x,r_N(u))} f(z) \, dz - \int\limits_{B(y,r_N(w))} f(z) \, dz \Big)^{N-2}.
\end{gather*}

Taking into account that $f(x) > 0$ and $f(y) > 0$ we establish that there exists
\begin{gather*}
    \widetilde{F}_{x,y}(u,w):=\lim_{\scriptscriptstyle{N \rightarrow \infty}} \widetilde{F}_{N,x,y}(u,w) = 1 - \lim_{\scriptscriptstyle{N \rightarrow \infty}} {\Big(1 - \frac{u\, f(x)}{\widetilde\gamma (N-1)} \Big)}^{N-2}  \\ - \lim_{\scriptscriptstyle{N \rightarrow \infty}} {\Big(1 - \frac{w\, f(y)}{\widetilde\gamma (N-1)} \Big)}^{N-2} + \lim_{\scriptscriptstyle{N \rightarrow \infty}} {\Big(1 - \frac{u\, f(x) + w \, f(y)}{\widetilde\gamma (N-1)} \Big)}^{N-2}
\end{gather*}
$$
    = 1 - e^{-\frac{f(x) u}{\widetilde\gamma}} - e^{-\frac{f(y) w}{\widetilde\gamma}} + e^{-\frac{f(x) u + f(y) w}{\widetilde\gamma}}\\
    = (1 - e^{-\frac{f(x) u}{\widetilde\gamma}}) (1 - e^{-\frac{f(y) w}{\widetilde\gamma}}) = F_x(u) F_y(w).
$$
Thus $\widetilde{F}_{x,y}(\cdot,\cdot)$ is the distribution function of a vector
$\eta _{x,y}:=(\eta_x,\eta_y)$
where $\eta_x\sim Exp(\frac{f(x)}{\widetilde{\gamma}})$,
$\eta_y\sim Exp(\frac{f(y)}{\widetilde{\gamma}})$
and the components of $\eta_{x,y}$ are independent. For $(x,y) \in A_1$, $i\neq j$ and $N\geq 3$, set
$$
\eta_{N,x}^{y,i,j}:= (N-1) V_d \widetilde{\gamma} \min{\left\{\min_{k\in\{1,\ldots,N\}\setminus \{i,j\}}\rho^d(x,X_k), \rho^d(x,y)\right\}}.
$$
Observe that the distribution function of $\eta_{N,x}^{y,i,j}$ is
$$
\p (\eta_{N,x}^{y,i,j} \leq u) = 1 - \ind\{ \rho(x,y) > r_N(u)\} (1-p_{N,x}(u))^{N-2}:=\widetilde{F}_{N,x}^y(u)
$$
for all $i,j\in \{1,\ldots,N\}$, $i\neq j$,
where $p_{N,x}(u)$ is defined by \eqref{p_int}.
Moreover $\widetilde{F}_{N,x,y}(\cdot,\cdot)$ is a distribution function of a random vector
$\eta_{N,x,y}^{i,j} : = (\eta_{N,x}^{y,i,j},\eta_{N,y}^{x,i,j})$, $(x,y) \in A_1$, $i,j\in \{1,\ldots,N\}$, $i\neq j$, $N \geq 3$.
Thus we have  shown that $\eta_{N,x,y}^{i,j}\stackrel{law}\rightarrow \eta_{x,y}$ as $N\to \infty$. Therefore
$$
\log\eta_{N,x}^{y,i,j}\log\eta_{N,y}^{x,i,j}\stackrel{law}\rightarrow \log\eta_x\log\eta_y,\;\;N\to \infty.
$$
Here we exclude a set of zero probability where random variables under consideration can be equal to zero.
Note that
\begin{gather*}
\e   \log \eta_{N,x}^{y,i,j} \, \log \eta_{N,y}^{x,i,j} = \int_{(0,\infty)} \int_{(0, \infty)} \log{u} \log{w} \, d \widetilde{F}_{N,x,y}(u, w) = \\ = \e  \big(\log{e^{\zeta_i(N)}} \log{e^{\zeta_j(N)}} | X_i = x, X_j = y \big) = \e  \big( \zeta_i(N) \zeta_j(N) | X_i = x, X_j = y \big).
\end{gather*}
Consequently, $\e \zeta_i(N)\zeta_j(N)=\e \zeta_1(N)\zeta_2(N)$ for all $i,j\in \{1,\ldots,N\}$, $i\neq j$. Hence in view of \eqref{t2: second: eq2} to prove that $\var H_n\to 0$ as $N\to \infty$
it is sufficient to show that $\cov (\zeta_1(N),\zeta_2(N))\to 0$, $N\to \infty$. Obviously,
$\e  \log{\eta_x} \log{\eta_y} = \e  \log{\eta_x} \, \e  \log{\eta_y} = \log{f(x)} \, \log{f(y)}$,
since $\eta_x$ and $\eta_y$ are independent.

For any fixed $M>0$, let us introduce $A_{1,M} := \big\{ (x,y) \in A_1: \rho(x,y) > M \big\}$. Our goal now is to prove, for any $M > 0$ and for all $(x,y)\in A_{1,M}$, that
\begin{align}\label{q1}
\begin{gathered}
    \e  \big( \zeta_1(N) \zeta_2(N) | X_1 = x, X_2 = y \big) \to \log{f(x)} \, \log{f(y)},\;\;N\to \infty.
\end{gathered}
\end{align}
First of all, we will establish the uniform integrability of a family $\{\log\eta_{N,x}^y\log\eta_{N,y}^x\}_{N \geq \widetilde{N}_0(x,y)}$ for such $(x,y) \in A_{1,M}$, where to simplify notation we write $\eta_{N,x}^y:= \eta_{N,x}^{y,1,2}$
and $\eta_{N,y}^x:= \eta_{N,y}^{x,1,2}$. The function $G(\cdot)$ is nondecreasing and convex. Thus,
for any $(x,y)\in A_{1,M}$ and $N\geq 3$,
\begin{align}\label{Step2: 1}
\begin{gathered}
    \e  G(|\log \eta_{N,x}^y \, \log \eta_{N,y}^x|) \leq \e G\left(\frac{(\log\eta_{N,x}^y)^2}{2} + \frac{(\log \eta_{N,y}^x)^2}{2}\right) \\ \leq \frac{1}{2} \left(\e G((\log \eta_{N,x}^y)^2) + \e G((\log \eta_{N,y}^x)^2)\right).
\end{gathered}
\end{align}
If $u \in (0,\frac{1}{e}]$, then one has $\widetilde{F}_{N,x}^y(u) = 1-(1-p_{N,x}(u))^{N-2}$ for any $N > \frac{1}{M^d e \widetilde\gamma V_d}+1$.
In similarity to \eqref{ineq_B} one can show that, if
$\rho(x,y) > M$ and $u \in (0, \frac{1}{e}]$, for any $N \geq \widetilde{N}_1$ ($\widetilde{N}_1$ does not depend on $x$ and $y$, but can depend on $M$, whence $\widetilde{N}_1 := \widetilde{N}_1(M)$),
    $\widetilde{F}_{N,x}^y(u) \leq \widetilde\gamma^{-\varepsilon_1}
    (M_f(x, R_1))^{\varepsilon_1} u^{\varepsilon_1}$.
Moreover, for all $u > 0$,
$$
1-\widetilde{F}_{N,x}^y(u) = \ind\{\rho(x,y)>r_N(u)\}(1-p_{N,x}(u))^{N-2} \leq (1-p_{N,x}(u))^{N-2}.
$$
Then analogously to \eqref{third: eq4} we derive that, for all $N \geq \widetilde{N}_2$ ($\widetilde{N}_2$ does not depend on $x$ and $y$) and $\varepsilon := \frac{\varepsilon_0 \varepsilon_2}{1+ \varepsilon_0}$ (we can assume as before that $\varepsilon \leq e$),
\begin{equation*}
    1 - \widetilde{F}_{N,x}^y(u)\leq \left(\frac{u}{2 \widetilde\gamma} m_f(x,R_2)\right)^{-\varepsilon}.
\end{equation*}
The same reasoning as was used at \textit{Step 1} of the proof of Theorem 2 leads, for $N\geq \max\{\widetilde{N}_1,\widetilde{N}_2\}$, to the inequality
\begin{align}\label{Step2: 2}
\begin{gathered}
    \e G(\log^2 \eta_{N,x}^y)  \leq \widetilde{C}_1 M_f^{\varepsilon_1}(x, R_1) +
    m_f^{-\varepsilon}(x, R_2)(\widetilde{C}_2  + \widetilde{C}_3(N)  + \widetilde{C}_4(N)  P(x))
\end{gathered}
\end{align}
where $\widetilde{C}_3(N) \to 0$ and $\widetilde{C}_4(N) \to 0$ as $N \to \infty$.
Then, in view of \eqref{Step2: 1} and \eqref{Step2: 2} for all $N \geq \max\{ \widetilde{N}_1, \widetilde{N}_2\}$,
\vspace{-0.5cm}
\begin{align}\label{Step2: 3}
\begin{gathered}
    \e  G(|\log \eta_{N,x}^y \, \log \eta_{N,y}^x|) \\ \leq \frac{\widetilde{C}_1}{2} (M_f^{\varepsilon_1}(x,R_1)+M_f^{\varepsilon_1}(y,R_1)) + \frac{\widetilde{C}_2}{2} (m_f^{-\varepsilon}(x,R_2)+m_f^{-\varepsilon}(y,R_2)) \\ + \frac{\widetilde{C}_3(N)}{2} (m_f^{-\varepsilon}(x,R_2) + m_f^{-\varepsilon}(y,R_2)) + \frac{\widetilde{C}_4(N)}{2} (m_f^{-\varepsilon}(x, R_2) P(x) + m_f^{-\varepsilon}(y, R_2) P(y)).
\end{gathered}
\end{align}
Consequently, in view of \eqref{Step2: 3}, for all $\kappa > 0$ and  $N \geq \widetilde{N}_0(x,y) = \max\{ \widetilde{N}_1, \widetilde{N}_2, \widetilde{N}_3(x,y)\}$,
\begin{gather*}
    \e  G(|\log \eta_{N,x}^y \, \log \eta_{N,y}^x|) \\ \leq \frac{\widetilde{C}_1}{2} (M_f^{\varepsilon_1}(x,R_1)+M_f^{\varepsilon_1}(y,R_1)) + \frac{\widetilde{C}_2}{2} (m_f^{-\varepsilon}(x,R_2)+m_f^{-\varepsilon}(y,R_2)) + \kappa := \widetilde{C}_0(x,y).
\end{gather*}
Hence, for any $(x,y) \in A_{1,M}$, a family $\{\log\eta_{N,x}^y \log\eta_{N,y}^x\}_{N \geq \widetilde{N}_0(x,y)}$ is uniformly integrable. Thus we come to \eqref{q1} for $(x,y)\in A_{1,M}$.

Let us define $T_N(x,y) := \e  \big( \zeta_1(N) \zeta_2(N) | X_1 = x, X_2 = y \big) = \e  \log \eta_{N,x}^y \, \log \eta_{N,y}^x$ where $(x,y)\in A_1$ and $N\geq 3$. The statement \eqref{q1} is equivalent to the following one: for any $(x,y)\in A_{1,M}$, $T_N(x,y) \to \log{f(x)} \log{f(y)}$ as $N \to \infty$. Take any $(x,y) \in A_1$. Then, for any fixed $M > 0$, we have proved that
\begin{equation}\label{Step2: asconv}
T_N(x,y) \ind\{\rho(x,y) > M\} \to \log{f(x)} \log{f(y)} \ind\{\rho(x,y) > M\}, \; N \to \infty.
\end{equation}
Note that
\begin{align}\label{Step2: 4}
\begin{gathered}
    G(|T_N(x,y)| \ind\{\rho(x,y) > M\}) \leq G(|T_N(x,y)|) = G(|\e  \log \eta_{N,x}^y \, \log \eta_{N,y}^x|) \\ \leq G(\e |\log \eta_{N,x}^y \, \log \eta_{N,y}^x| ) \leq \e G(|\log \eta_{N,x}^y \, \log \eta_{N,y}^x|).
\end{gathered}
\end{align}
Due to \eqref{Step2: 3} and \eqref{Step2: 4} one can conclude that, for all $N \geq \{\widetilde{N}_1, \widetilde{N}_2 \}$,
\begin{gather*}
\int_{\mathbb{R}^d} \int_{\mathbb{R}^d} G(|T_N(x,y)| \ind\{\rho(x,y) > M\}) f(x) f(y) \, dx \, dy \\ \leq \widetilde{C}_1 \int_{\mathbb{R}^{d}} M_f^{\varepsilon_1}(x, R_1) f(x) \, dx + \widetilde{C}_2 \int_{\mathbb{R}^{d}} m_f^{-\varepsilon}(x, R_2) f(x) \, dx \\ + \widetilde{C}_3(N) \int_{\mathbb{R}^d} m_f^{-\varepsilon}(x, R_2) f(x) \, dx +  \widetilde{C}_4(N) \int_{\mathbb{R}^d} m_f^{-\varepsilon}(x, R_2) P(x) f(x) \, dx.
\end{gather*}
Therefore, for any $\kappa > 0$, there exists $\widetilde{N}_4(\kappa) >0$ such that
\begin{gather*}
\int_{\mathbb{R}^d} \int_{\mathbb{R}^d} G(|T_N(x,y)| \ind\{\rho(x,y) > M\}) f(x) f(y) \, dx \, dy \\ \leq \widetilde{C}_1 Q_f(\varepsilon_1, R_1) + \widetilde{C}_2 T_f(\varepsilon, R_2) + \kappa < \infty
\end{gather*}
for $N \geq \widetilde{N}_0 := \max\{\widetilde{N}_1(M), \widetilde{N}_2, \widetilde{N}_4(\kappa)\}$, since $\varepsilon = \frac{\varepsilon_0 \varepsilon_2}{1+\varepsilon_0} < \varepsilon_2$.
Hence, for $(x,y)\in A_1$, a family
$\big\{ T_N(x,y) \ind\{ \rho(x,y) > M\}\big\}_{N \geq \widetilde{N}_0}$ is uniformly integrable w.r.t. $\p_{\xi}\otimes \p_{\xi}$.
Consequently, in view of \eqref{Step2: asconv}, for each $M>0$,
\begin{equation}\label{final_conv}
\int_{\rho(x,y) > M}\!\!\!\! T_N(x,y) f(x) f(y) \, dx \, dy \to \int_{\rho(x,y) > M}\!\!\!\! \log{f(x)} \log{f(y)} f(x) f(y) \, dx \, dy, \;\; N \to \infty.
\end{equation}
Now we consider the case $\rho(x,y) \leq M$. For each $M < 1$, one has
$$
\p\Big(\rho(X_1,X_2) \leq M\Big) = \p\Big(G(|\log\rho(X_1,X_2)|) \geq G(|\log{M}|)\Big)
\leq \frac{\e G(|\log\rho(X_1,X_2)|)}{G(|\log{M}|)} \to 0
$$
as $M\to 0$, since $\e G(|\log\rho(X_1,X_2)|) < \infty$ due to the condition $K_{f,2}(\varepsilon_0)<\infty$.

In view of the elementary inequality $t \leq G(t) + 1$ for all $t \geq 0$ we can write
\begin{gather*}
    \big|\e \zeta_1(N) \zeta_2(N) \ind\{\rho(X_1,X_2) \leq M\} \big| \\ = \left|\int_{\mathbb{R}^d} \int_{\mathbb{R}^d} \e (\zeta_1(N) \zeta_2(N) \ind\{\rho(X_1,X_2) \leq M\} | X_1 = x, X_2 = y) f(x) f(y) \, dx \, dy \right|
    \\
    \leq \int_{\mathbb{R}^d} \int_{\mathbb{R}^d} \e \big|\log \eta_{N,x}^y \log \eta_{N,y}^x\big| \ind\{\rho(x,y) \leq M\}  f(x) f(y) \, dx \, dy \\ \leq \int_{\rho(x,y) \leq M} \e G(|\log \eta_{N,x}^y \log \eta_{N,y}^x|) f(x) f(y) \, dx \, dy +
    \p\Big(\rho(X_1,X_2) \leq M\Big).
\end{gather*}

Hence by \eqref{Step2: 3}, for all $N \geq \max\{\widetilde{N}_1, \widetilde{N}_2\}$, one has
\begin{gather*}
    \int_{\rho(x,y) \leq M} \e G(|\log \eta_{N,x}^y \log \eta_{N,y}^x|) f(x) f(y) \, dx \, dy \\ \leq \widetilde{C}_1 \; \e M_f^{\varepsilon_1}(X_1,R_1) \ind\{\rho(X_1,X_2) \leq M\}  + \widetilde{C}_2 \; \e m_f^{-\varepsilon}(X_1,R_2) \ind\{\rho(X_1,X_2) \leq M\} \\ + \widetilde{C}_3(N) \; \e m_f^{-\varepsilon}(X_1,R_2) \ind\{\rho(X_1,X_2) \leq M\} + \widetilde{C}_4(N) \; \e m_f^{-\varepsilon}(X_1, R_2) P(X_1) \ind\{\rho(X_1,X_2) \leq M\} \\ = \widetilde{C}_1 \widetilde{C}_5(M) + \widetilde{C}_2 \widetilde{C}_6(M) + \widetilde{C}_3(N) \widetilde{C}_7(M) + \widetilde{C}_4(N) \widetilde{C}_8(M)
\end{gather*}
where $\widetilde{C}_k(M) \to 0$ as $M \to 0$ for $k = 5,6,7,8$. Here we take into account that the Lebesgue integral is absolutely continuous function ($\e M_f^{\varepsilon_1}(X_1,R_1) < \infty$, $\e m_f^{-\varepsilon}(X_1,R_2) < \infty$, $\e m_f^{-\varepsilon}(X_1, R_2) P(X_1) < \infty$) and
$\p\left(\rho(X_1,X_2) \leq M\right)\to 0$ as $M\to 0$.

Hence, for any $\kappa > 0$, one can find $M_1=M_1(\kappa) > 0$ such that, for all $M \in (0,M_1]$ and $N \geq \max\{\widetilde{N}_1, \widetilde{N}_2, \widetilde{N}_5(\kappa)\}$,
\begin{gather*}
    \left| \int_{\rho(x,y) \leq M} T_N(x,y) f(x) f(y) \, dx \, dy \right| = \big|\e \zeta_1(N) \zeta_2(N) \ind\{\rho(X_1,X_2) \leq M\} \big| < \frac{\kappa}{3}.
\end{gather*}
Also  (from absolute continuity of the Lebesgue integral)
one can find $M_2=M_2(\kappa)>0$ such that, for all $M\in (0,M_2]$,
$$
\left| \int_{\rho(x,y) \leq M}\log{f(x)} \log{f(y)} f(x) f(y) \, dx  \, dy \right| < \frac{\kappa}{3}
$$
Take $M= \min\{M_1,M_2\}$. Due to \eqref{final_conv} one can find $\widetilde{N}_6(M, \kappa)$ such that for all $N \geq \widetilde{N}_6(M, \kappa)$ the following inequality holds
\begin{gather*}
    \left| \int_{\rho(x,y) > M} T_N(x,y) f(x) f(y) \, dx \, dy - \int_{\rho(x,y) > M}
\log{f(x)} \log{f(y)} f(x) f(y) \, dx \, dy \right| < \frac{\kappa}{3}.
\end{gather*}
So, for any $\kappa > 0$ one can find $M(\kappa) >0$ and $\widetilde{N}_0(\kappa) :=\max\{\widetilde{N}_1, \widetilde{N}_2, \widetilde{N}_5(\kappa), \widetilde{N}_6(M, \kappa)\}$ such that for all $N \geq \widetilde{N}_0(\kappa)$:
\begin{gather*}
  \bigg|\int_{\mathbb{R}^{d}} \int_{\mathbb{R}^{d}} T_N(x,y) f(x) f(y) \, dx \, dy -
\int_{\mathbb{R}^{d}} \int_{\mathbb{R}^{d}} \log{f(x)} \log{f(y)} f(x) f(y) \, dx \, dy \bigg|
< \kappa.
\end{gather*}
Therefore $\e \zeta_1(N)\zeta_2(N)\to H^2$, $N\to \infty$, and consequently
$\cov(\zeta_1(N),\zeta_2(N))\to 0$ as $N\to \infty$ because $\lim_{N\to \infty}\e \zeta_1(N)
= \lim_{N\to \infty}\e \zeta_2(N) =H$.
The proof is complete. $\square$

\vskip0.2cm
{\it Proof of Corollary \ref{cor3}}.
In view of the proof of Corollary \ref{cor2}, it remains to demonstrate that condition \eqref{p} for some $p > 2$ implies the finiteness of $K_{f,2}(\varepsilon_0)$ for some $\varepsilon_0 > 0$.

Assume that condition \eqref{p} is satisfied for some $p>2$. In Section 2 we used the elementary inequality:
for any $z\geq 1$ and $\delta >0$
$\log z \leq \frac{1}{\delta}z^{\delta}$.
Hence, for any $t \geq 0$ and $\delta >0$, $G(t^2)\leq \frac{1}{\delta}t^{2(1+\delta)}$.
By the Lyapunov inequality, for any $\delta>0$ and $\varepsilon_0 >0$,
$$
\left(\int_{\mathbb{R}^d}|\log\rho(x,y)|^{2(1+\delta)}
f(y)\,dy\right)^{1+\varepsilon_0}\leq \int_{\mathbb{R}^d}|\log\rho(x,y)|^{2(1+\delta)(1+\varepsilon_0)}
f(y)\,dy.
$$
For a given $p>2$ take $\delta:=\frac{p-2}{p+2}$ and $\varepsilon_0:=\frac{p-2}{4}$ then
$2(1+\delta)(1+\varepsilon_0)=p$. Thus
\begin{gather*}
K_{f,2}(\varepsilon_0) \leq \left(\frac{1}{\delta}\right)^{1+\varepsilon_0} \int_{\mathbb{R}^d} \left(\int_{\mathbb{R}^d}|\log\rho(x,y)|^{2(1+\delta)}
f(y)\,dy\right)^{1+\varepsilon_0} f(x) \, dx \\ \leq
\left(\frac{1}{\delta}\right)^{1+\varepsilon_0} \int_{\mathbb{R}^d}\int_{\mathbb{R}^d}|\log\rho(x,y)|^{p}
f(x) f(y)\,dx\,dy < \infty.
\end{gather*}
We conclude that condition \eqref{p} for $p > 2$  implies the required finiteness: $K_{f,2}(\varepsilon_0)<\infty$ for $\varepsilon_0=\frac{p-2}{4}$. \;$\square$

\vskip0.2cm
{\it Proof of Corollary \ref{cor4}}.
Let $\xi \sim N(\nu, \Sigma)$, where $\nu \in \mathbb{R}^d$ and $\Sigma$ is a (strictly) positive-definite $d\times d$ matrix. Thus $\xi$ has a density
\begin{equation}\label{normal_density}
    f(x) = \frac{1}{(2 \pi)^{\frac{d}{2}} (\det\Sigma)^{\frac{1}{2}}} e^{-\frac{1}{2}(\Sigma^{-1}  (x-\nu),x - \nu)},\;\;x\in \mathbb{R}^d,
\end{equation}
where $(\cdot,\cdot)$ stands for a scalar product in $\mathbb{R}^d$ related to the
Euclidean norm $\|\cdot\|$.

First of all we prove that for such $f$ relation \eqref{reg} holds.
Denote by $\{\lambda_1,\cdots,\lambda_d\}$ the set of all eigenvalues of a matrix $\Sigma$.
Then $\lambda_{min} := \min_{i=1,\dots,d}\lambda_i >0$
since $\lambda_i > 0$
for all $i \in \{ 1,\ldots,d\}$.
Thus $\Sigma = C^{\top}\Lambda C$ for a diagonal matrix
$\Lambda := \text{diag}(\lambda_1,\ldots,\lambda_d)$ and an orthogonal matrix $C$. Hence,  $\Sigma^{-\frac{1}{2}} := C^T \Lambda^{-\frac{1}{2}} C$, where $\Lambda^{-\frac{1}{2}} = \text{diag}(\lambda_1^{-1/2}, \ldots, \lambda_d^{-1/2})$. Now one can rewrite \eqref{normal_density} as follows
\begin{equation}\label{normal_density2}
    f(x) = \frac{1}{(2 \pi)^{\frac{d}{2}} (\det\Sigma)^{\frac{1}{2}}} e^{-\frac{1}{2} \|\Sigma^{-\frac{1}{2}} (x - \nu)\|^2},\;\;x\in \mathbb{R}^d.
\end{equation}
Note also that, for all $v, \, w \in \mathbb{R}^d$, the relation
$\|v\|^2 = \|w\|^2 + \|v-w\|^2 + 2(w,v-w)$ holds. Therefore
\begin{align}\label{eeeee}
\begin{gathered}
\|\Sigma^{-\frac{1}{2}} (z - \nu)\|^2 = \|\Sigma^{-\frac{1}{2}} (x - \nu)\|^2 + \|\Sigma^{-\frac{1}{2}} (z-x)\|^2 + 2 \big( \Sigma^{-\frac{1}{2}} (x-\nu), \Sigma^{-\frac{1}{2}} (z-x) \big) \\ \leq \|\Sigma^{-\frac{1}{2}} (x - \nu)\|^2 + \frac{1}{\lambda_{min}} \|z-x\|^2 + 2 (\Sigma^{-1} (x-\nu), z-x )
\end{gathered}
\end{align}
where the last inequality is true because $\|\Sigma^{-\frac{1}{2}} (z-x)\|^2
\leq \frac{1}{\lambda_{min}} \| z-x\|^2$.
According to \eqref{normal_density2} and \eqref{eeeee} we have for all $z,\, x \in \mathbb{R}^d$
\begin{equation}\label{gauss}
    f(z) \geq
    f(x) \,  e^{-\frac{1}{2 \lambda_{min}} \| z-x\|^2} e^{(\Sigma^{-1} (\nu-x), \,z-x)}.
\end{equation}
Let us fix arbitrary $x \in \mathbb{R}^d$, $R > 0$ and take any $r$ such that  $0 < r \leq R$. By \eqref{gauss} we get
\begin{align}\label{1_final}
\begin{gathered}
    \int_{B(x,r)} f(z) \, dz \geq f(x) \, \int_{B(x,r)} e^{-\frac{1}{2 \lambda_{min}} \| z-x\|^2} e^{(\Sigma^{-1} (\nu-x), \,z-x)} \, dz \\ \geq
    f(x) \, e^{-\frac{r^2}{2 \lambda_{min}}} \, \int_{B(x,r)} e^{(\Sigma^{-1} (\nu-x), \,z-x)} \, dz \geq f(x) \, e^{-\frac{R^2}{2 \lambda_{min}}} \, \int_{B(0,r)} e^{(\Sigma^{-1} (\nu-x) \,y)} \, dy
\end{gathered}
\end{align}
because $\| z-x\|^2 \leq r^2$ for $z \in B(x,r)$.
Simple inequality $e^t \geq 1+t$, $t \in \mathbb{R}$, leads to formula
\begin{align}\label{2_final}
\begin{gathered}
    \int_{B(0,r)} e^{(\Sigma^{-1} (\nu-x), \,y)} \, dy \geq \int_{B(0,r)} \big(1+(\Sigma^{-1} (\nu-x), \,y)\big) \, dy \\ = r^d V_d + \int_{B(0,r)} (\Sigma^{-1} (\nu-x), \,y)\, dy = r^d V_d
\end{gathered}
\end{align}
since by the Fubini theorem $\int_{B(0,r)} (v, y) \, dy = 0$ for any fixed $v \in \mathbb{R}^d$ (as $\int_{(-a,a)} u \, du = 0$ for any $a > 0$).
\vspace{-0.1cm}
In view of \eqref{1_final} and \eqref{2_final} we have $\int_{B(x,r)} f(z) \, dz \geq f(x) \, e^{-\frac{R^2}{2 \lambda_{min}}} \, r^d V_d$. Consequently, for any $0 < r \leq R$ and any $x \in \mathbb{R}^d$, we come to \eqref{reg} with $c = e^{-\frac{R^2}{2 \lambda_{min}}}$.
Note that \eqref{densi} holds as, for $\varepsilon\in (0,1)$,
\begin{equation*}
\int_{\mathbb{R}^d} f^{1-\varepsilon}(x) \, dx  =  \left(\frac{(2\pi)^{\varepsilon}}{1-\varepsilon}\right)^{\frac{d}{2}}<\infty.
\end{equation*}
Then for an arbitrary Gaussian random vector, the integral in \eqref{p3} is finite for all $\varepsilon_2 \in (0,1)$ and all $R_2 > 0$.

Clearly, (B) is valid with $M= (2\pi)^{-\frac{d}{2}}(\det{\Sigma})^{-\frac{1}{2}}$. So, \eqref{p2} is finite for any $\varepsilon_1 > 0$ and $R_1 > 0$.

Now we prove that, for each $p>1$, condition \eqref{p} is satisfied, thus condition (A) is satisfied as well. In other words we show that
$\e |\log{\norm{\xi_1 - \xi_2}}|^p <\infty$ where $\xi_1$ and $\xi_2$ are independent copies of $\xi$. One can write
$$
\e |\log{\norm{\xi_1 - \xi_2}}|^p
=\e |\log{\norm{\xi_1 - \xi_2}}|^p \mathds{1}\{ \norm{\xi_1 - \xi_2} \leq 1\} + \e |\log{\norm{\xi_1 - \xi_2}}|^p \mathds{1}\{ \norm{\xi_1 - \xi_2} > 1\}.
$$
For each $x\in \mathbb{R}^d$, one has
$$\e (|\log{\norm{\xi_1 - \xi_2}}|^p \mathds{1}\{ \norm{\xi_1 - \xi_2} \leq 1\}|\xi_2 = x) = \e |\log{\norm{\xi_1 - x}}|^p \mathds{1}\{ \norm{\xi_1 - x} \leq 1\}
$$
$$
= \e (-\log{\norm{\eta_x}})^p \mathds{1}\{ \norm{\eta_x} \leq 1\}$$
where $\eta_x = \xi - x$, $\eta_x \sim N(\nu-x, \Sigma)$.
Taking the spherical coordinates we infer that
$$
\e (-\log{\norm{\eta_x}})^p \mathds{1}\{ \norm{\eta_x} \leq 1\}
\leq \frac{1}{(2 \pi)^{\frac{d}{2}} (\det\Sigma)^{\frac{1}{2}}} \int_{\norm{y} \leq 1} (-\log{\norm{y}})^p \, dy
$$
$$
= \frac{\Gamma(p+1)}{2^{d/2-1}d^{p+1} (\det\Sigma)^{\frac{1}{2}} \Gamma(\frac{d}{2})} := C(p,d).
$$
Hence $\e |\log{\norm{\xi_1 - \xi_2}}|^p \mathds{1}\{ \norm{\xi_1 - \xi_2} \leq 1\} \leq  C(p,d) < \infty$.

For $t>1$,  $|\log{t}| = \log{t} < t$ and, for $p>1$ and $u,w\in \mathbb{R}$, one has
$|u+w|^p\leq 2^{p-1}(|u|^p+|w|^p)$. Thus
$$
\e |\log{\norm{\xi_1 - \xi_2}}|^p \mathds{1}\{ \norm{\xi_1 - \xi_2} > 1\} \leq \e \norm{\xi_1 - \xi_2}^p \leq  \e (\norm{\xi_1-\nu} + \norm{\xi_2-\nu})^p
$$
$$
\leq 2^{p-1} (\e \norm{\xi_1-\nu}^p + \e \norm{\xi_2-\nu}^p) = 2^p \e \norm{\xi-\nu}^p< \infty.
$$
as one can easily verify that
\vspace{-0.2cm}
$$
\e \norm{\xi-\nu}^p \leq \frac{2^{\frac{p}{2}}}{(\det\Sigma)^{\frac{1}{2}}} \frac{\Gamma(\frac{p+d}{2})}
{\Gamma(\frac{d}{2})} (\lambda_{max})^{\frac{p+d}{2}} <\infty.
\vspace{-0.1cm}
$$
Consequently, the finiteness of \eqref{p1} and \eqref{p12} are established for any $\varepsilon_0 > 0$. We do not use here the
explicit formula for $H$.
The proof of Corollary \ref{cor4} is complete. \;$\square$
\vskip0.2cm
{\bf Acknowledgements} The work is supported by the Russian Science Foundation under grant 14-21-00162 and performed at the Steklov Mathematical Institute of Russian Academy of Sciences.
The authors are grateful to Professor E.Spodarev for drawing their attention to the entropy estimation problems.

\section{Appendix}

{\it Proof of Lemma \ref{l6a}}.
Let us fix $x_0\in \mathbb{R}^d$, $r_0>0$. For $B\in \mathcal{B}(\mathbb{R}^d)$, set $p(B):=\int_B f(y)dy$.
The function $r^{-d}$ is continuous on $(0,\infty)$.  Thus it is sufficient to verify that
$|p(B(x,r))-p(B(x_0,r_0))|$ is small whenever $(x,r)$ is close to $(x_0,r_0)$. One has
\begin{align}\label{eq1y}
\begin{gathered}
|p(B(x,r))-p(B(x_0,r_0))|\leq |p(B(x,r))-p(B(x,r_0))|+|p(B(x,r_0))-p(B(x_0,r_0))|\\
\leq p(B(x,r)\triangle B(x,r_0)) + p(B(x,r_0)\triangle B(x_0,r_0)).
\end{gathered}
\end{align}
Since $f\in L^1(\mathbb{R}^d)$, for any $\varepsilon >0$,
there exists $\delta =\delta(\varepsilon)>0$ such that $p(B) < \varepsilon$ if $\mu(B)<\delta$.
Clearly, $\mu(B(x,r)\triangle B(x,r_0)) =  V_d|r^d- r_0^d|$ and
$$
B(x,r_0)\triangle B(x_0,r_0)\subset
(B(x_0,r_0+\|x-x_0\|)\setminus B(x_0,r_0)) \cup (B(x,r_0+\|x-x_0\|)\setminus B(x,r_0)).
$$
Thus
$\mu(B(x,r_0)\triangle B(x_0,r_0))\leq 2V_d((r_0 + \|x-x_0\|)^d - r_0^d)$. Taking into account
\eqref{eq1y} we come to the statement of Lemma. $\square$

\vskip0.2cm
{\it Proof of Lemma \ref{lemma1}}.
1) The Lyapunov inequality (i.e. the Jensen inequality for a power function) yields that,
for any $\varepsilon \in (0,\varepsilon_0]$,
$K_f(\varepsilon) \leq K_f(\varepsilon_0)^{\frac{1+\varepsilon}{1+\varepsilon_0}} <\infty$.

2) Assume that $Q_f(\varepsilon_1, R_1) < \infty$. Consider $Q_f(\varepsilon_1, R)$ where $R>0$.
If $0<R\leq R_1$ then, for each $x\in \mathbb{R}^d$, according to the definition of $M_f$ one has
$M_f(x,R)\leq M_f(x,R_1)$. Hence $Q_f(\varepsilon_1, R)\leq Q_f(\varepsilon_1, R_1)<\infty$.
Let now $R>R_1$. One has
$$
M_f(x,R)
\leq \max\left\{M_f(x,R_1), \sup_{R_1<r\leq R}\frac{\int_{B(x,R_1)}f(x)dx + \int_{B(x,r)\setminus B(x,R_1)}f(x)dx}{|B(x,r)|}\right\}
$$
$$
\leq \max\left\{M_f(x,R_1), M_f(x,R_1)+\frac{1}{|B(x,R_1)|}\right\} = M_f(x,R_1)+\frac{1}{|B(x,R_1)|}.
$$
Therefore
\begin{gather*}
Q_f(\varepsilon_1, R) = \int_{\mathbb{R}^d} (M_f(x,R))^{\varepsilon_1} f(x) \, dx \leq \int_{\mathbb{R}^d} \left(M_f(x,R_1)+\frac{1}{R_1^d V_d} \right)^{\varepsilon_1} f(x) \, dx \\ \leq  \max\{1,2^{\varepsilon_1-1}\}\left(Q_f(\varepsilon_1, R_1) + \left(\frac{1}{R_1^d V_d}\right)^{\varepsilon_1}\right)< \infty.
\end{gather*}

Assume now that $Q_f(\varepsilon_1,R)<\infty$ for some $\varepsilon_1>0$ and $R>0$. Then, for
any $\varepsilon \in (0,\varepsilon_1]$,
the Lyapunov inequality yields $Q_f(\varepsilon,R)\leq (Q_f(\varepsilon_1,R))^{\frac{\varepsilon}{\varepsilon_1}}<\infty$.

3) Now we turn to $T_f(\varepsilon,R)$, $\varepsilon\in (0,\varepsilon_2]$, $R>0$. Let $T_f(\varepsilon_2, R_2) < \infty$ and
take $0<R\leq R_2$. Then, for each $x$, according to the definition of $m_f$ we get
$0\leq m_f(x,R_2) \leq m_f(x,R)$. Hence $T_f(\varepsilon_2,R)\leq T_f(\varepsilon_2,R_2)<\infty$.
Consider $R>R_2$.  According to Lemma \ref{l6a}, for each $x\in \mathbb{R}^d$ and every $a>0$, the function
$I_f(x,r)$
is continuous in $r$ on $(0,a]$. If $x\in S(f)\cap \Lambda(f)$ then
there exists $\lim_{r\to 0+}I_f(x,r)=f(x)$. For such $x$, set $I_f(x,0):=f(x)$. Further in the proof we consider an arbitrary fixed $x\in S(f)\cap \Lambda(f)$. Then $I_f(x,\cdot)$ is continuous on any interval $[0,a]$. Thus one can find $\widetilde{R}_2$ in $[0,R_2]$ such that
$m_f(x,R_2)= I_f(x,\widetilde{R}_2)$ and there exists $R_0$ in $[0,R]$ such that
$m_f(x,R)=I_f(x,R_0)$. If $R_0 \leq R_2$ then $m_f(x,R)=m_f(x,R_2)$ (since $m_f(x,R)\leq m_f(x,R_2)$ for $R>R_2$ and $m_f(x,R)=I_f(x,R_0) \geq m_f(x,R_2)$ as $R_0 \in [0,R_2]$).
Assume that
$R_0 \in (R_2,R]$. Obviously $R_0>0$ as $R_2>0$. One has
\begin{gather*}
m_f(x,R)= I_f(x,R_0)=\frac{\int_{B(x,R_2)}f(y)dy+\int_{B(x,R_0)\setminus B(x,R_2)}f(y)dy}{|B(x,R_0)|} \\
\geq \frac{\int_{B(x,R_2)}f(y)dy}{|B(x,R_0)|}= \frac{|B(x,R_2)|}{|B(x,R_0)|}I_f(x,R_2)\geq \frac{|B(x,R_2)|}{|B(x,R_0)|}m_f(x,R_2) \\
= \left(\frac{R_2}{R_0}\right)^d m_f(x, R_2) \geq \left(\frac{R_2}{R}\right)^d m_f(x, R_2).
\end{gather*}
Thus
for all cases ($R_0 \in [0, R_2]$ and $R_0\in (R_2,R]$) one has $m_f(x,R) \geq \left(\frac{R_2}{R}\right)^d m_f(x,R_2)$ as $R_2 < R$. Taking into account Remark 1 and the relation
$\mu(S(f)\setminus (S(f)\cap \Lambda(f)))=0$ we come to the inequality
\vspace{-0.5cm}
$$
T_f(\varepsilon_2,R)\leq \left(\frac{R}{R_2}\right)^{\varepsilon_2 d}T_f(\varepsilon_2,R_2)<\infty.
$$
Assume now that $T_f(\varepsilon_2,R)<\infty$ for some $\varepsilon_2>0$ and $R>0$. Then, for
any $\varepsilon \in (0,\varepsilon_2]$,
the Lyapunov inequality yields $T_f(\varepsilon,R)\leq (T_f(\varepsilon_2,R))^{\frac{\varepsilon}{\varepsilon_2}}<\infty$. The proof is complete. $\square$

\vskip0.3cm

{\it Proof of Lemma \ref{lemma_G}}. We start with relation 1).
Note that if a function $g$ is measurable and bounded on a finite interval $(a,b]$ then
$\int_{(a,b]}g(x)\nu(dx)$ is finite where $\nu$ is a finite measure on the Borel subsets of $(a,b]$. Thus, for each $a\in (0,\frac{1}{e}]$, using the integration by parts formula (see, e.g., \cite{Shiryaev}, p. 245) we get
\begin{align}\label{lGcor1}
\begin{gathered}
\int_{(a, \frac{1}{e}]} F(u) \frac{\log{(-\log{u})} + 1}{u} \, du =
\int_{(a, \frac{1}{e}]} F(u)d(-(-\log u)\log(-\log u)) \\
= (-\log a) \log (-\log a)F(a) + \int_{(a, \frac{1}{e}]}(-\log u)\log(-\log u)d F(u).
\end{gathered}
\end{align}
Assume now that
$\int_{(0, \frac{1}{e}]} (-\log{u}) \log{(-\log{u})}\,dF(u)<\infty$.
Then by the monotone convergence theorem
\begin{equation}\label{c1}
\lim_{a\to 0+}\int_{(0, a]}(-\log u)\log(-\log u)d F(u)=0.
\end{equation}
Clearly, the following nonnegative integral admits an estimate
\begin{gather*}
\int_{(0, a]} (-\log{u}) \log{(-\log{u})} \, dF(u) \geq (-\log{a}) \log{(-\log{a})} \int_{(0, a]} dF(u)  \\
= (-\log{a}) \log{(-\log{a})} (F(a) - F(0)) =
(-\log{a}) \log{(-\log{a})} F(a) \geq 0.
\end{gather*}
Therefore \eqref{c1} implies that
\begin{equation}\label{c2}
(-\log{a}) \log{(-\log{a})} F(a) \to 0,\;\;a\to 0+.
\end{equation}
Taking $a\to 0+$ in \eqref{lGcor1} we come by the monotone convergence theorem to  relation 1) of the statements of our Lemma.

Suppose now that
\begin{equation}\label{c3}
\int_{(0, \frac{1}{e}]} F(u) \frac{\log{(-\log{u})} + 1}{u} \, du <\infty.
\end{equation}
In view of \eqref{lGcor1}, \eqref{c3} and monotone convergence theorem we have
$$
\lim_{b\to 0+}\int_{(0, b]}F(u)d(-(-\log u)\log(-\log u))\,du=0.
$$
For any $c\in (0,b)$ we obtain the inequalities
$$
\int_{(0, b]}F(u)d(-(-\log u)\log(-\log u))
\geq \int_{(c, b]}F(u)d(-(-\log u)\log(-\log u))
$$
$$
= F(b)(-(-\log b)\log (-\log b)) + F(c)(-\log c)\log (-\log c)
+ \int_{(c,b]}(-\log u)\log(-\log u)\,dF(u)
$$
$$
\geq F(b)(\log b)\log (-\log b) - F(c)(\log c)\log (-\log c)
+ (-\log b)\log(-\log b)(F(b)-F(c))
$$
$$
= F(c)(-\log c) \log(-\log c)\left(1 - \frac{\log b \log(-\log b)}{\log c \log(-\log c)}\right).
$$
Let $c=b^2$ ($b\leq\frac{1}{e} < 1$). Then
$$
1 - \frac{\log b \log(-\log b)}{\log c \log(-\log c)}= \frac{1}{2} + \frac{\log2}{2\log(-2\log b)}\geq \frac{1}{2}
$$
for all positive $b$ small enough. Hence we can claim that
$F(b^2)(-\log(b^2))\log(-\log(b^2))\to 0$ as $b\to 0$. Therefore we come to \eqref{c2} taking $a=b^2$.
Consequently we get \eqref{c2}. Then \eqref{lGcor1} yields relation 1).

If one of (nonnegative) integrals appearing in 1) is infinite and other one is finite
we come to the contradiction. Hence 1) is established.
In a similar way one can prove that relation 2) is valid. Therefore, we omit further details.
$\square$

\vskip0.3cm
{\it Proof of Lemma \ref{l1}}.
Take $x\in S(f)\cap \Lambda(f)$ and $R>0$. Assume that $m_f(x,R)=0$.
According to Lemma \ref{l6a}, there exists $\widetilde{R}\in [0,R]$
($\widetilde{R}=\widetilde{R}(x,R)$) such that $m_f(x,R)= I_f(x,\widetilde{R})$ (recall that  $I_f(x, 0) := \lim_{r \rightarrow 0+} I_f(x,r) = f(x)$ for all $x \in \Lambda(f)$ by continuity).
If $\widetilde{R}=0$ then $m_f(x,r)=f(x)>0$ as $x\in S(f)\cap \Lambda(f)$.
Hence we have to consider $\widetilde{R} \in (0,R]$.
If $I_f(x,\widetilde{R})=0$ then $\int_{B(x,r)}f(y)dy =0$ for any $0<r\leq \widetilde{R}$.
Thus \eqref{2} ensures that $f(x)=0$. However, $x\in S(f)\cap \Lambda(f)$.
Consequently, for all $x \in S(f) \cap \Lambda(f)$ we have $m_f(x,R) > 0$. Hence, $S(f)\cap \Lambda(f)\subset D_f(R):=\{x\in S(f): m_f(x,R)>0\}$. It remains to note
that $S(f)\setminus \Lambda(f)\subset \mathbb{R}^d \setminus \Lambda(f)$ and $\mu(\mathbb{R}^d \setminus \Lambda(f))=0$. Therefore $\mu(S(f) \setminus D_f(R)) = 0$. $\square$

\vskip0.3cm
{\it Proof of Lemma \ref{l4}}.
We verify that, for a given $c>0$ and any $t>0$ ($t:=\rho^d(x,y)$), there exist $D,F\geq 0$ such that
\begin{equation}\label{m1}
G(|\log (ct)|)\leq D G(|\log t|)+F.
\end{equation}
Instead of \eqref{m1} one can prove that, for any $v\in \mathbb{R}$ ($v:= \log t$) and
$c_1:=\log c \in \mathbb{R}$,
$
G(|v+c_1|)\leq D G(|v|)+F.
$
Obviously $|v+c_1|\leq |v|+|c_1|$. Since $G$ is non decreasing function it is enough to verify that
$G(|v|+|c_1|)\leq D G(|v|)+F$.
Thus we prove that ($s:=|v|$, $c_2 := |c_1|$)
\begin{equation}\label{e1}
G(s+c_2)\leq DG(s)+F,\;\;s\geq 0, c_2\geq 0.
\end{equation}
Let $0\leq s \leq 2$. Then $G(s+c_2)\leq G(2+c_2)$ and, for $s\in [0,2]$, the inequality \eqref{e1} holds with any $D= 0$ and
$F=G(2+c_2)$.
Let $s>2$. Then $G(s)=s\log s$, $G(s+c_2)= (s+c_2)\log(s+c_2)$ and
$$
\frac{(s+c_2)\log(s+c_2)}{s\log s} = \left(1+\frac{c_2}{s}\right)\frac{\log(s+c_2)}{\log s}
\leq \left(1+\frac{c_2}{2}\right)\frac{\log(s+c_2)}{\log s}.
$$
One can find $T_0>2$ such that $\frac{\log(s+c_2)}{\log s}\leq 2$ for $s > T_0$ as
$\lim_{s\to \infty}\frac{\log(s+c_2)}{\log s}=1$.
Obviously, $
\frac{\log(s+c_2)}{\log s}\leq\frac{\log (T_0+c_2)}{\log 2}
$
for $2<s \leq T_0$.
Consequently we come to \eqref{e1} with $D=(1+\frac{c_2}{2})\max\{2,\frac{\log (T_0+c_2)}{\log 2}\}$ and $F=G(2+c_2)$.

Now we establish that, for some $D_0,F_0\geq 0$ and all $t>0$,
$
G(|\log t^d|)\leq D_0 G(|\log t|)+F_0.
$
Note that $G(|\log t^d|)= G(d|\log t|)$. Thus it is sufficient to show that,
for a given $c\geq 0$ and any $v \geq 0$ ($v:=|\log t|$) and $c\geq 0$, one has
$G(cv)\leq D_0G(v)+F_0$. This is verified in a similar way to \eqref{e1} proof.
Thus, for some $a,b\geq 0$ (actually $a := D D_0$, $b := D F_0 + F$) and all $x,y \in \mathbb{R}^d$, $x\neq y$, we come to the inequality $G(|\log\widetilde{\rho}(x,y)|)\leq a G(|\log\rho(x,y)|) +b$ and consequently
\eqref{v0} is valid (since $\int_{\mathbb{R}^d}f(y)dy=1$). $\square$

\end{document}